\documentclass[twocolumn]{autart}    

\usepackage{lineno,hyperref}								
\usepackage{indentfirst,dsfont}								
\usepackage{color,parskip}											
\usepackage{graphicx}										

\usepackage{mathtools,balance}
\usepackage{caption}
\usepackage[section]{placeins}
\usepackage{amssymb}
\usepackage{indentfirst}
\usepackage{psfrag}
\usepackage{tabularx}
\usepackage{savesym}
\savesymbol{AND}
\usepackage{algorithm,cite}
\usepackage{algorithmic,balance}

\usepackage{wasysym,subcaption,multirow}
\usepackage{graphicx}

\graphicspath{{figures/}}

\newtheorem{assump}{Assumption}
\newtheorem{lemma}{Lemma}
\newtheorem{theo}{Theorem}

\DeclareMathOperator*{\minimise}{minimise}
\DeclareMathOperator*{\argmin}{argmin}

\newcommand{\red}[1]{\textcolor{red}{#1}}

\begin{document}

\begin{frontmatter}
	
	\title{Subgradient averaging for multi-agent optimisation with different constraint sets} 
	
	\thanks[footnoteinfo]{L. Romao is supported by the Coordination for the Improvement of Higher Education Personnel (CAPES) - Brazil. The work of K. Margellos and A. Papachristodoulou has been supported by EPSRC UK under grants $ \mathrm{EP}/\mathrm{P}03277\mathrm{X}/1  $ and $ \mathrm{EP}/\mathrm{M}002454/1 $, respectively. Giuseppe Notarstefano is supported by the European Research Council (ERC) under the European Union’s Horizon $ 2020 $ research and innovation programme (grant $ 638992  $- OPT4SMART).}
	
	\author[Oxford]{Licio Romao}\ead{licio.romao@eng.ox.ac.uk},    
	\author[Oxford]{Kostas Margellos}\ead{kostas.margellos@eng.ox.ac.uk},               
	\author[Bologna]{Giuseppe Notarstefano}\ead{giuseppe.notarstefano@unibo.it.}  
	\author[Oxford]{Antonis Papachristodoulou}\ead{antonis@eng.ox.ac.uk}  
	
	\address[Oxford]{Department of Engineering Science, University of Oxford, Parks Road, Oxford OX1 3PJ, UK.}  
	\address[Bologna]{Department of Electrical, Electronic, and Information Engineering G. Marconi at Alma Mater Studiorum Università di Bologna.}             


	\begin{abstract}                          
			
		We consider a multi-agent setting with agents exchanging information over a possibly time-varying network, aiming at minimising a separable objective function subject to constraints. To achieve this objective we propose a novel subgradient averaging algorithm that allows for non-differentiable objective functions and different constraint sets per agent. Allowing different constraints per agent simultaneously with a time-varying communication network constitutes a distinctive feature of our approach, extending existing results on distributed subgradient methods. To highlight the necessity of dealing with a different constraint set within a distributed optimisation context, we analyse a problem instance where an existing algorithm does not exhibit a convergent behaviour if adapted to account for different constraint sets. For our proposed iterative scheme we show asymptotic convergence of the iterates to a minimum of the underlying optimisation problem for step sizes of the form $  \frac{\eta}{k+1} $, $ \eta > 0 $. We also analyse this scheme under a step size choice of $ \frac{\eta}{\sqrt{k+1}} $, $ \eta > 0 $, and establish a convergence rate of $ \mathcal{O}(\frac{\ln k}{\sqrt{k}}) $ in objective value. To demonstrate the efficacy of the proposed method, we investigate a robust regression problem and an $ \ell_2 $ regression problem with regularisation.
	\end{abstract}
	
\end{frontmatter}

\section{Introduction}
\label{sec:introduction}

Distributed optimisation deals with multiple agents interacting over a network and has found numerous applications in different domains, such as wireless sensor networks~\cite{MG:12,BMG:14}, robotics~\cite{MBCF:07}, and power systems~\cite{BCCZ:15}, due to its ability to parallelize computation and prevent agents from sharing information considered as private. Typically, distributed algorithms are based on an iterative process in which agents maintain some estimate about the decision vector in an optimisation context, exchange this information with neighbouring agents according to an underlying communication protocol/network, and update their estimate on the basis of the received information.

Despite the intense research activity in this area, only a few algorithms can simultaneously deal with time-varying networks, non-differentiable objective functions and account for the presence of constraints~\cite{ZM:12,NO:15,XK:17,MFGP:18,LWY:19}, features that are often treated separately in the literature. Several of the commonly employed methods are based on a projected subgradient or a proximal step and their analysis consists of selecting the step size underlying these algorithms, establishing a convergence rate analysis, and quantifying practical convergence for (near-)real time applications.

In this paper, we study a class of optimisation problems that involves a separable objective function, while the feasible set can be decomposed as an intersection of different compact convex sets. A centralised version of this class of problems has been studied under a stochastic setting in \cite{Bia:2016,PN:18}. Distributed algorithms for this class have been proposed in~\cite{JKJJ:08,NO:09,NOP:10,ZM:12,LRS:16,LN:13,MFGP:18,MA:19}. References \cite{JKJJ:08,NO:09,NOP:10} rely on \cite{TBA:86,BT:89} to propose a distributed strategy based on projected sub-gradient methods. These results consist of an averaging step followed by a local sub-gradient projection update.  In \cite{MFGP:18} a distributed scheme based on a proximal update is proposed, thus extending \cite{JKJJ:08,NOP:10} to the case where different local constraint sets and an arbitrarily time-varying network are considered. The authors in \cite{ZM:12} provide asymptotic convergence for a primal-dual algorithm that allows coupling between agents' local estimates. We discuss additional related results in Section \ref{sec:Related_work}, after the proposed algorithm is presented and some notation introduced.

We motivate our approach by constructing an example showing that extending available
	algorithms to the case of different constraint sets might
	not exhibit a convergent behaviour for all problem instances.  Hence, a direct adaptation of existing schemes is not always possible when dealing with different constraint sets. Notice also that   distributed algorithms developed for the unconstrained case cannot be trivially adapted to our setting, as lifting the constraints in the objective (e.g., via characteristic functions) would violate boundedness of the subgradient, a typical requirement for such algorithms \cite{DAW:12,NOP:10,NO:15,MFGP:18}.

The main contribution of this paper is the introduction and the characterization of the convergence rate for a new subgradient averaging algorithm. The proposed scheme allows us to account for time-varying networks, non-differentiable objective functions and different constraint sets per agent as in~\cite{MFGP:18}, while achieving faster practical convergence as it is based on subgradient averaging as in~\cite{JKJJ:08,DAW:12,MA:19}. Note that allowing simultaneously for different constraint sets per agent and time-varying communication network by means of a subgradient averaging scheme is a distinct feature of the algorithm in this paper.  Preliminary results related to this paper appeared in \cite{RMNP:19b}, where several proofs have been omitted. Moreover, the construction of Section \ref{subsec:counter_ex} that motivates the analysis of algorithms with different constraint sets is novel, and offers insight on the limitations of existing algorithms. We also provide detailed numerical examples, not included in the conference version.

The paper is organised as follows. In Section~\ref{sec:pro_statement} we present the problem statement, the network communication structure, and the main assumptions adopted in this paper, followed by a numerical construction that motivates the algorithm of this paper. In Section~\ref{sec:distributed_method} we present the proposed scheme and the main convergence results, namely, asymptotic convergence in iterates and a convergence rate as far as the optimal value is concerned. Section \ref{sec:Related_work} provides detailed discussion and comparison of our scheme with other results in the literature. In Section~\ref{sec:numerical_examples} we study the robust linear regression problem and $ \ell_2 $ regression with regularisation to demonstrate the main algorithmic features of our scheme and to compare our strategy against existing methods. Finally, some concluding remarks and future research directions are provided in Section~\ref{sec:conclusion}. To ease the reader all proofs have been deferred to the Appendix (Section~\ref{sec:appendix}).

\emph{Notation:} We denote by $ \mathbb{R} $ the set of real numbers and $ \mathbb{N} $ the set of natural numbers (excluding zero). The symbol $ \mathbb{R}^n $ stands for the Cartesian product $ \mathbb{R} \times \ldots \times \mathbb{R} $ with $ n $ terms. A sequence of elements in $ \mathbb{R}^n $ is denoted by $ (x(k))_{k \in \mathbb{N}} $. For any set $ X \subset \mathbb{R}^n$, we denote its interior, relative interior and convex hull by $ \mathrm{int}(X) $, $ \mathrm{ri}(X) $, and $ \mathrm{conv}(X) $, respectively. We also denote by $ f(X) $ as the image of the set $ X $ over a function $ f $. The subdifferential of $ f $ at a point $ x \in \mathrm{dom} f $ is denoted by $ \partial f(x) $. For any point $ x \in \mathbb{R}^n $, $ \| x \|_2 $ stands for the Euclidean norm of $ x $ and $ \| x \|_1 $ for the $ \ell_1 $ norm of $ x \in \mathbb{R}^n $, which are reduced to $ |x| $ if $ x $ is scalar.
\section{Problem statement and a motivating example}
\label{sec:pro_statement}
\subsection{Problem set-up and network communication}
\label{subsec:pro_setup}
Consider the optimisation problem
\begin{equation}
\begin{aligned}
\minimise_x & \quad f(x) = \sum_{i = 1}^m f_i(x) \\ 
\mbox{subject to} & \quad x \in \cap_{i = 1}^m X_i, 
\end{aligned}
\label{eq:main_pro}
\end{equation}
where $ x \in \mathbb{R}^n $ is the vector of decision variables, and $ f_i: \mathbb{R}^n \rightarrow \mathbb{R} $ and $ X_i \subset \mathbb{R}^n $ constitute the local objective function and constraint set, respectively, for agent $ i $, $ i = 1, \ldots, m $. We suppose that each agent $ i $ possesses as private information the pair $ (f_i,X_i) $ and maintains a local estimate $ x_i $ of the common decision vector $ x $. 

The goal is for all agents to agree on the local variables, that is, $ x_i = x^\star $, for all $ i = 1, \ldots, m $, where $ x^\star $ is an optimiser of~\eqref{eq:main_pro}, i.e., a feasible point such that $ f(x^\star) \leq f(x)  $ for all $ x \in \cap_{i = 1}^m X_i $. We impose the following assumption.

\begin{assump} We assume that:
	\begin{itemize}
		\item[$ i) $] For all $ i = 1, \ldots, m $, the function $ f_i $ is convex.
		
		\item[$ ii) $] The set $ X_i \subset \mathbb{R}^n $ is compact and convex for all $ i = 1, \ldots, m $, and $ \cap_{i = 1}^m X_i  $ has a non-empty interior.
		
		\item[$ iii) $] The subgradient of the function $ f(x) $ is bounded on $ \cup_{i = 1}^m X_i $, that is, $  
		L = \max_{\substack{\xi \in \partial f(x),\\ x \in \cup_{i = 1}^m X_j}} \| \xi \|_2 < \infty.
		$
	\end{itemize}
	\label{assump:Convexity_Compt}
\end{assump}

Assumption~\ref{assump:Convexity_Compt} imposes standard restriction for constrained non-smooth optimisation. Item $ ii) $ implies informally that $ \cup_{i = 1}^m X_i $ has volume in $ \mathbb{R}^n $, i.e., that the affine hull of $ \cup_{i = 1} X_i $ has dimension $ n $. Moreover, the compactness assumption of item $ ii) $ guarantees that the optimal set of problem~\eqref{eq:main_pro} is non-empty. Item $ iii) $ is an assumption that is needed to prove convergence of sub-gradient methods applied to problem~\eqref{eq:main_pro}. Under item $ iii) $, the sub-gradient of the function $ f $ can be evaluated at points that belong to $ \cup_{i = 1}^m X_i $. We provide in Appendix~\ref{app:suff_bounded_subgrad} a technical condition on the domain of the functions $ f_i $ that is sufficient to guarantee that Assumption~\ref{assump:Convexity_Compt}, item $ iii) $, holds. An important consequence of Assumption~\ref{assump:Convexity_Compt} is given in the following lemma.

\begin{lemma}
	Under Assumption~\ref{assump:Convexity_Compt}, we have that:
	\begin{itemize}
		\item[$ i) $] The set $ \mathrm{conv}(\cup_{i = 1}^m X_i) $ is compact.
		\item[$ ii) $] The function $ f $ is Lipschitz continuous over $ \cap_{i = 1}^m X_i $, i.e., the following inequality hods 
		\[ 
		| f(x) - f(y) | \leq L \| x - y \|_2, \quad \forall~x,y~\in~\cap_{i = 1}^m X_i,
		\]
		where $ L $ is the constant defined in Assumption~\ref{assump:Convexity_Compt}.
	\end{itemize}
	\label{lemma:bounded_subgradient}
\end{lemma}

Typical choices of functions that satisfy Assumption~\ref{assump:Convexity_Compt} are piecewise-linear functions, quadratic convex functions and the logistic regression function.

In this paper, we aim to solve problem~\eqref{eq:main_pro} through a network of agents that use only the available local information, namely, the pair $ (f_i, X_i) $ and the current estimate for the optimal solution, $ x_i(k) $, $ i = 1, \ldots, m $, maintained by agent $ i $ at a given instance $ k $. We will show how $ x_i(k), i = 1, \ldots, m,$ can be constructed and updated in Section \ref{sec:distributed_method}, with $ k $ playing the role of iteration index. To this end, we now characterise the underlying communication network.  Let $  \mathcal{G}(k) = (\mathcal{N}, \mathcal{E}(k))  $ be an undirected graph, where $ \mathcal{N} = \{ 1, \ldots, m \} $ is the number of agents  and $ \mathcal{E}(k) \subset \mathcal{N} \times \mathcal{N} $ is the set of edges at iteration $ k $, that is, only if node $ (j,i) \in \mathcal{E}(k) $ then node $ j $ sends information to node $ i $ at iteration $ k $. We associate the time-varying matrix $ A(k) $ to the edge set $ \mathcal{E}(k) $, with $ [A(k)]_j^i > 0 $ only if $ (j,i) \in \mathcal{E}(k) $ at time $ k $. As the graph is undirected, the matrix $ A(k) $ can be chosen to be symmetric. We also define the graph $ \mathcal{G}_{\infty} = (\mathcal{N}, \mathcal{E}_{\infty}) $, in which $ (j,i) \in \mathcal{E}_{\infty} $ if agent $ j $ communicates with agent $ i $ infinitely often. We impose the following assumption on the matrix $ A(k) $.

\begin{assump} We assume that:
	\begin{itemize}
		\item[$ i) $] The graph $ (\mathcal{N}, \mathcal{E}_{\infty}) $ is strongly connected. Moreover, there exits a uniform upper bound on the communication time for all $ (j,i) \in \mathcal{E}_{\infty} $.
		\item[$ ii) $] There exists $ \eta \in (0,1) $ such that for all $ k \in \mathbb{N} $ and for all $ i,j = 1, \ldots, m $, $ [A(k)]_i^i \geq \eta $, and if $ [A(k)]_j^i > 0 $ then we have that  $ [A(k)]_j^i $ $\geq \eta$.
		\item[$ iii) $]  Matrix $ A(k) $ is doubly stochastic.
	\end{itemize}
	\label{assump:Strongly_graph}
\end{assump}

These are standard requirements in the distributed optimisation literature. We refer the reader to~\cite{NO:09b,NOP:10,DAW:12,MFGP:18} for more details.

\subsection{Dealing with different constraint sets}
\label{subsec:counter_ex}
In this section, we highlight the necessity of developing a new algorithmic scheme to deal with the presence of a different constraint sets per agent. To this end, consider the iterative scheme\footnote{It should be noted that $ z_i, i = 1, \ldots, m,  $ in~\eqref{eq:aveg_step_mod_Duchi} should not be confused with that of Step~2 in Algorithm~\ref{algo:distributed_algorithm} presented in the sequel; we use the same symbol to match the notation in~\cite{DAW:12} and ease the reader.}
\begin{subequations}
		\label{eq:mod_Duchi}
	\begin{align}
	z_i(k+1) &= \sum_{j = 1}^{m} [A]_j^i z_j(k) + g_i(k) \label{eq:aveg_step_mod_Duchi} \\
	x_i(k+1) &= \argmin_{\xi \in X_i} z_i(k+1)^T \xi  + \frac{1}{c(k)} \| \xi \|_2^2 \label{eq:x_update_mod_Duchi},
	\end{align}
\end{subequations}

\noindent which consists of a modified version of the algorithm considered in~\cite{DAW:12}, adapted to account for different constraint sets in each agent's local optimisation problem. In the setting of the previous section, notice that matrix $ A $ in~\eqref{eq:aveg_step_mod_Duchi} corresponds to a time-invariant network $ \mathcal{G}(k) = (\mathcal{N}, \mathcal{E}) $, for all $ k \in \mathbb{N} $. Assumption~\ref{assump:Strongly_graph} is satisfied if the graph $ (\mathcal{N}, \mathcal{E}) $ is strongly connected and matrix $ A $ is doubly-stochastic.

Observe that~\eqref{eq:aveg_step_mod_Duchi} constitutes a subgradient update step, with neighbouring local variables $ z_j(k) $ being ``mixed'' according to the matrix $ A$ and added to $ g_i(k) \in \partial f_i(x_i(k)) $, i.e., a subgradient of $ f_i $ evaluated at $ x_i(k) $, $ i = 1, \ldots, m $. Step~\eqref{eq:x_update_mod_Duchi} is an optimisation program with the objective function being the sum (weighted via $ c(k) $) of 
\begin{equation*}
	z_i(k+1)^T \xi \text{: linear ``proxy'' of }  f_i, 
	\end{equation*} and a regularization term $ \| \xi \|_2^2 $. To comply with~\cite{DAW:12}, we set $ c(k) = \frac{1}{\sqrt{k+1}} $. Recall that the algorithm in~\cite{DAW:12} involves the same constraint set in the update rule of~\eqref{eq:x_update_mod_Duchi}, that is $ X_i = X $ for all $ i = 1, \ldots, m $, and possesses a guaranteed convergence rate of $ \mathcal{O}(\frac{\ln k}{\sqrt{k}}) $ for the running averages of the iterates $ x_i(k) $; here, we introduce a different set $ X_i $ per agent and show that this (natural) modification may lead to erroneous results. 

\begin{figure}[!htb]
	\centering
	\includegraphics[width=0.85\linewidth]{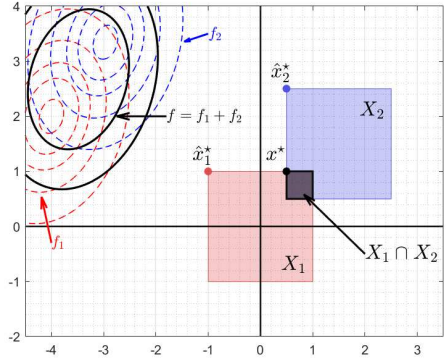}
	\caption{Geometric representation of problem instance encoded by~\eqref{eq:problem_matrices}. The red ellipsoids (dashed lines) correspond to the level curves of $ f_1 $, the blue ellipsoids (double-dashed lines) represent the function $ f_2 $, while the black (solid lines) ellipsoids to the ones of $ f = f_1 + f_2 $. The shaded red box illustrates the constraint set $ X_1 $, while the shaded blue box illustrates  $ X_2 $. Vectors $ \hat{x}^\star_1  = [-1,1]^T$  and $ \hat{x}^\star_2 = [0.5,2.5]^T$ are the optimal solutions of $ f_1(x) $ and $ f_2(x) $ under the constraints $ X_1 $ and $ X_2 $, respectively. The global optimal solution of $ f = f_1 + f_2 $ with matrices given by~\eqref{eq:problem_matrices} subject to $ x \in X_1 \cap X_2 $ is denoted by~$ x^\star $. This construction shows that $ \hat{x}_1^\star$ and $ \hat{x}_2^\star $ constitute fixed-points of~\eqref{eq:mod_Duchi} thus preventing the iteration from reaching $ x^\star $ if initialised at those points.} 
	\label{fig:Counter_ex}
\end{figure}

Consider a two-agent instance of~\eqref{eq:main_pro}, i.e., $ m = 2 $ with $ x \in \mathbb{R}^2 $, $ f_i = x^T Q x + q_i^T x + r_i, $ for $ i = 1,2  $ and
\begin{align}
Q = \begin{bmatrix}
1.2 & 0.4 \\
0.4 & 1.8
\end{bmatrix}, \quad q_1 &= \begin{bmatrix}
8 \\ -4
\end{bmatrix}, \quad q_2 = \begin{bmatrix}
2.93 \\ -11.46
\end{bmatrix}, \nonumber \\ 
r_1 &= 20, \quad r_2 = 25.
\label{eq:problem_matrices}
\end{align}
The local constraint sets are given by $ X_1 = [-1,1] \times [-1,1] $ and $ X_2 = [0.5,2.5] \times [0.5,2.5] $. The feasible set $ X_1 \cap X_2 $ is the box $ [0.5,1] \times [0.5,1] $. Figure~\ref{fig:Counter_ex} depicts the level curves of the quadratic functions $ f_1(x) $ (dashed-red lines), $ f_2 $ (double-dashed lines), and $ f = f_1 + f_2 $ (solid-black lines). The red and blue boxes represent the sets $ X_1 $ and $ X_2 $ respectively, with the feasible set, $ X_1 \cap X_2 $, being also indicated in the figure in black.

By inspection the optimal solution of $ f_1 $ under the constraint $ x \in X_1 $ is $\hat{x}^\star_1 = \left[ 
-1, 1
\right]^T$. Similarly, the optimal solution for $f_2$ under $x \in X_2 $ is $\hat{x}^\star_2 = \left[ 
0.5, 2.5
\right]^T$. We then have the following proposition. 

\begin{prop}
	Let $ (z_i(k))_{k \in \mathbb{N}}, (x_i(k))_{k \in \mathbb{N}} $, $ i = 1, 2 $, be the sequences generated by algorithm~\eqref{eq:mod_Duchi} when applied to problem~\eqref{eq:problem_matrices} with initial conditions $ x_i(0) = \hat{x}_i^\star $, $ i = 1, 2 $, and with $ A = \frac{1}{2} \mathbf{1} \mathbf{1}^T $ and $ c(k) = \frac{1}{\sqrt{k+1}} $. We have that
	\[ 
	x_1(k) = \hat{x}^\star_1, \quad x_2(k) = \hat{x}^\star_2, \quad \forall k \in \mathbb{N}.
	\]
	\label{prop:Counter_ex}
\end{prop}

Proposition~\ref{prop:Counter_ex} shows that $ \hat{x}^\star_1 $ and $ \hat{x}^\star_2 $ constitute fixed points of~\eqref{eq:mod_Duchi}, hence the iteration cannot reach $ x^\star $ if initialised from these points. This highlights the necessity of devising a new algorithm to deal with the presence of a different constraint set per agent.
 \section{Distributed Methodology}
 \label{sec:distributed_method}
\subsection{Proposed algorithm}
\label{subsec:prop_algo}
The main steps of the proposed scheme are summarized in Algorithm~\ref{algo:distributed_algorithm}. We initialise each agents' local variable with an arbitrary $ x_i(0) \in X_i $, $ i = 1, \ldots, m $; such points are not required to belong to $ \cap_{i = 1}^m X_i $.

At iteration $ k $, agent $ i $ receives $ x_j $ from the neighbouring agents and averages them through $ A(k) $, which captures the communication network, to obtain $ z_i(k) $. Recall that we denote the element of the $ j $-th row and $ i $-th column of matrix $ A(k) $ by $ [A(k)]_j^i $. Agent $ i $ then calculates a subgradient, $ g_i $, of its own objective function evaluated at $ z_i(k) $ and broadcasts this information back to its neighbours. In the sequel, agent $ i $ averages the received $ g_j(z_j(k)) $ in order to compose a proxy for a subgradient of $ f(x) $, namely, $ d_i(k) $. Finally, agents minimise a linear proxy $ d_i(k)^T \xi  $ of $ f(\xi) $ plus a regularization term weighted by $ \frac{1}{c(k)} $. An alternative interpretation of this last computation is that agents update their local estimates by performing a subgradient step with step size $ c(k) $ and projecting $ z_i(k) - c(k) d_i(k)$ onto their local set. Indeed, this local update can be rewritten as
\[ 
x_i(k+1) = \mathcal{P}_{X_i} [z_i(k) - c(k) d_i(k) ]
\]
where $ \mathcal{P}_{X_i} [\cdot] $ denotes projection onto the set $ X_i $. 

\begin{algorithm}[!htbp]
	\caption{Proposed distributed algorithm}
	\small	
	\begin{algorithmic}[] 
		\REQUIRE: $ x_i(0), \quad  i = 1, \ldots, m $ 
		\vspace{0.25cm}
		\STATE \textbf{For} $ i = 1, \ldots, m $, \textbf{repeat until convergence}
		\vspace{0.1cm}
	\end{algorithmic}
\begin{algorithmic}[1]
		\STATE Compute $ z_i(k) = \sum_{j = 1}^{m} [A(k)]_j^i  x_j(k),$
		\vspace{0.25cm}
		\STATE Pick $ g_i(z_i(k)) \in \partial f_i(z_i(k)),$
		\vspace{0.25cm}
		\STATE Compute $ d_i(k) = \sum_{j = 1}^{m} [A(k)]_j^i g_j(z_j(k)),$
		\vspace{0.25cm}
		\STATE Compute $ x_i (k + 1) = \argmin_{\xi \in X_i} d_i(k)^T \xi + \frac{1}{2 c(k)} \| z_i(k) - \xi \|_2^2, $ 
		\vspace{0.25cm}
		\STATE Set $ k \leftarrow k +1  $
		\end{algorithmic}	
	\begin{algorithmic}
		\vspace{0.1cm}
		\STATE \textbf{~~end}
	\end{algorithmic}
	\normalsize
	\label{algo:distributed_algorithm}
\end{algorithm}

\subsection{Algorithm Analysis}
\label{sec:algo_analysis}

\subsubsection{Convergence in iterates}
\label{subsec:convergence_iterates}

In this subsection, we impose the following assumption on the step size $ c(k) $.

\begin{assump}
	\label{assump:time_varying_step_size_1_t}
	Let $ (c(k))_{k \in \mathbb{N}} $ be the sequence adopted in Algorithm~\ref{algo:distributed_algorithm}. We require that:
	\begin{itemize}
		\item[$ i) $] $ c(k) $ is non-negative and non-increasing;
		\item[$ ii) $] $ \sum_{k = 1}^\infty c(k) = \infty$ and $ \sum_{k = 1}^\infty c(k)^2 < \infty $.
	\end{itemize}
\end{assump}

A sequence satisfying Assumption~\ref{assump:time_varying_step_size_1_t} is $ c(k) = \frac{\eta}{k + 1} $, for $ \eta > 0 $. 

\begin{theo}
	Let $ (x_i(k))_{k \in \mathbb{N}} $ be the sequences generated by Algorithm~\ref{algo:distributed_algorithm}, for all $ i = 1, \ldots, m $. Under Assumptions~\ref{assump:Convexity_Compt}-~\ref{assump:time_varying_step_size_1_t}, we have that for some minimizer $ x^\star  $ of~\eqref{eq:main_pro},
	\[
	\lim_{k \rightarrow \infty} \| x_i(k) - x^\star \|_2 = 0, \quad \forall~i = 1, \ldots, m. 
	\]
	\label{theo:convergence_algo_1}
\end{theo}
The proof of Theorem~\ref{theo:convergence_algo_1}, as well as of Theorem~\ref{theo:convergence_rate} presented in the sequel, is based on some auxiliary technical results presented in Appendix~\ref{app:Auxiliary_results}.

Theorem~\ref{theo:convergence_algo_1} extends the result in~\cite{MFGP:18} by allowing an agent to communicate subgradient information to neighbouring agents, a feature that, as illustrated in Section~\ref{sec:numerical_examples}, can speed up practical convergence.

\subsubsection{Convergence in objective value and convergence rate}
\label{subsec:conv_value_conv_rate}

Throughout this section, we impose the following assumption on the step size $ c(k) $.

\begin{assump}
	The sequence $ (c(k))_{k \in \mathbb{N}} $ used in Algorithm~\eqref{algo:distributed_algorithm} is $ c(k) = \frac{\eta}{\sqrt{k+1}} $, for some $ \eta > 0 $.
	\label{assump:time_varying_step_size_1_sqrt}
\end{assump}

Our convergence rate results build on the following related sequence generated by Algorithm~\ref{algo:distributed_algorithm}, 
\begin{equation}
\hat{x}_i(k+1) = \frac{c(k+1) x_i(k+1) + S(k)\hat{x}_i(k)}{S(k+1)}, 
\label{eq:running_average}
\end{equation}
where $ S(k) = \sum_{r = 1}^k c(r) $, and $ (x_i(k))_{k \in \mathbb{N}} $, for all $ i = 1, \ldots, m $, are the sequences generated by Algorithm~\ref{algo:distributed_algorithm}, with initial condition $ \hat{x}_i(0) = x_i(0) $. By rewriting~\eqref{eq:running_average} as $ 
\hat{x}_i(k) = \frac{1}{S(k)} \sum_{r = 1}^k c(r) x_i(r),
$ we can interpret \eqref{eq:running_average} as a convex combination of past iterates.
	\renewcommand{\arraystretch}{1.4}
	\begin{table*}[!htb]
		\scriptsize
		\caption{Summary of distributed schemes for smooth and non-smooth optimisation.}
		\centering 
		\begin{tabular}{c|c|cc|c || c|cc|c} 
			\hline\hline 
			&\multicolumn{4}{c||}{Smooth + Constant step-size} & \multicolumn{4}{c}{Non-smooth + Diminishing step-size} \\ [0.8ex] \cline{2-9} 
			&\multicolumn{2}{c}{Common sets} & \multicolumn{2}{c||}{Different sets} & \multicolumn{2}{c}{Common sets} & \multicolumn{2}{c}{Different sets} \\  [0.8ex]  \cline{2-9} 
			\multirow{2}{*}{} &Convex	& \multicolumn{1}{c||}{Strongly} & Convex	&\multicolumn{1}{c||}{Strongly}&	Convex	&\multicolumn{1}{c||}{Strongly}&	Convex& Strongly \\
			&	& \multicolumn{1}{c||}{Convex} & 	& \multicolumn{1}{c||}{Convex}&		& \multicolumn{1}{c||}{Convex}&	&  \multicolumn{1}{c}{Convex}  \\
			\hline 
			No (sub)grad. avg.  & \cite{NO:09,JMX:12,YLY:16} & \multicolumn{1}{c||}{\cite{YLY:16}}    & \cite{LCF:16} & \multicolumn{1}{c||}{-}   & \cite{NO:15,SS:18}  & \multicolumn{1}{c||}{\cite{TLR:12,LQX:17}} & \cite{ZM:12,LN:13,LRS:16,MFGP:18} & - \\ [0.8ex] \hline
			(Sub)grad. avg. & \cite{SLWY:15,QL:17,SS:18,VZCPS:16} & \multicolumn{1}{c||}{\cite{SLWY:15,QL:17,SS:18}}   & -  & - &  \cite{DAW:12,SS:18,XK:17,LWY:19} & \multicolumn{1}{c||}{-} & our work, \cite{MA:19} & - \\ [0.8ex] 
			\hline
			\hline 
		\end{tabular}
		\label{table:survey}
		\normalsize
\end{table*}
\begin{theo}
	Consider the running average defined in~\eqref{eq:running_average}. Under Assumptions~\ref{assump:Convexity_Compt},~\ref{assump:Strongly_graph}, and~\ref{assump:time_varying_step_size_1_sqrt}, we have that:
	\begin{itemize}
		\item[$ i) $] For all $ i,j = 1, \ldots, m $, the sequence $ (\| \hat{x}_i(k) - \hat{x}_j(k) \|)_{k \in \mathbb{N}} $ converges to zero at a rate $ \mathcal{O}(\frac{\ln k}{\sqrt{k}}) $.
		\item[$ ii) $] All accumulation points of the sequence $ (\hat{x}_i(k))_{k \in \mathbb{N}} $ are feasible.
		\item[$ iii) $] There exist $ B_1, B_2 > 0 $ such that 
		\begin{equation}
		\begin{aligned}
		\left|\sum_{i = 1}^m  f_i(\hat{x}_i(k)) - f(x^\star)\right| &\leq B_1 \frac{1}{\sqrt{k}} + B_2 \frac{\ln k}{\sqrt{k}}.
		\end{aligned} 
		\label{eq:ineq_x_convergence_rate}
		\end{equation}
	\end{itemize}
	\label{theo:convergence_rate}
\end{theo}
Note that Theorem~\ref{theo:convergence_rate} asserts convergence of the function value along the running average $ \hat{x}_i(k) $, i.e., all limit point of $ (\hat{x}_i(k))_{k \in \mathbb{N}} $ are optimal, however, the iterates might exhibit an oscillatory behaviour. For the exact expression of $ B_1 $ and $ B_2 $, we refer the reader to  Appendix~\ref{app:Proof_theo_2}.
The absolute value in Theorem~\ref{theo:convergence_rate} is due to the fact that $\hat{x}_i(k)$ may not be necessarily feasible; however, item $ ii) $ in Theorem~\ref{theo:convergence_rate} implies that all accumulation points of $ (\hat{x}_i(k))_{k \in \mathbb{N}} $, $ i = 1, \ldots, m $, are feasible. Item $ i) $ states the rate at which consensus is achieved for the sequences $ (\hat{x}_i(k))_{k \in \mathbb{N}} $, $  i = 1, \ldots, m $. Similar rates can be obtained with more general choices for the step size, e.g., $ c(k) = \frac{1}{k^a} $, for $ a \in [0.5,1). $

It should be noted that the result of Theorem~\ref{theo:convergence_rate} further extends the work presented in~\cite{MFGP:18} not only by allowing agents to communicate their (sub-) gradients, but by also unveiling how to (non trivially) adapt the proof line in that paper to come up with convergence results that recover traditional rates for distributed subgradient methods. This is the first convergence rate result under the scenario considered in this paper. 

\section{Comparison with related algorithms}

\label{sec:Related_work}

In this section we provide a detailed comparison of the proposed algorithm with other results in the literature. To this end, note that in \cite{JKJJ:08} a similar distributed sub-gradient scheme is mentioned, but no analysis of such a scheme is presented. References \cite{LRS:16,LN:13} characterize the convergence rate of a sub-gradient algorithm under different constraint sets per agent that does not possess subgradient averaging. References \cite{MFGP:18,ZM:12} show asymptotic convergence of distributed algorithms with different constraint sets and time-varying communication network. Hence, by combining (sub)-gradient averaging and providing an analysis that yields convergence rates under time-varying communication networks and different constraint sets per agent, the results in this paper are distinct from all the above literature.

A closely related algorithm to the one presented here is the one in \cite{MA:19}. This provides convergence rates assuming a regularity condition on the local sets (weaker than compactness) and requiring the network to be row-stochastic; however, it does not analyse the case where the communication network is time-varying. This requires different analysis arguments, thus complementing the results in \cite{MA:19}, extending them to allow for time-varying networks. Moreover, the example of Section \ref{subsec:counter_ex} highlights the need for developing a different analysis when agents possess different constraints sets.

Although only marginally related to the results of this paper, it is worth mentioning distributed algorithms that deal with similar optimization problems~\cite{SLWY:15,QL:17,SS:18}. Paper \cite{SS:18} proposes an algorithm whose convergence is valid for non-convex objectives and directed communication network, while \cite{SLWY:15,QL:17} use a constant step size to establish linear convergence rates for strongly convex functions. Moreover, distributed algorithms based on proximal methods with constant step sizes have been proposed in \cite{CO:12}. In this setting, the objective function is assumed to be differentiable to obtain convergence to the optimal solution of problem \eqref{eq:main_pro}, and the size of the allowable step-size is upper bounded by a quantity related to the Lipschitz constant of the objective function. Unlike these results, we allow for non-differentiable objective functions.
	
To better position this paper within the recent literature, we summarise the main distributed algorithms that are amenable to smooth and non-smooth constrained optimisation in Table~\ref{table:survey}. We highlight both scenarios of common and different local constraint sets, which are indicated in the table by common sets and different sets, respectively. In this brief summary, we restrict our attention to algorithms that use constant step size for smooth optimisation, and to those that use diminishing step sizes for the non-smooth case. We also present a categorization of these schemes between those that have results for general convex functions and strongly convex functions. In row entitled ``No (sub)grad. avg.'', we include distributed algorithms based on projected (sub)gradient, proximal minimisation, and primal-dual update that do not leverage on averaging first-order information from neighbouring agents. In contrast, row ``(Sub)grad. avg.'' includes algorithms that exploit (sub) gradient averaging. Among the few papers that are suitable for different local sets, this is the first result to establish a convergence rate that matches that of the common local sets case, and simultaneously allows agents to use first-order information of their neighbours under time-varying communication networks, thus speeding up practical convergence.
\section{Numerical Examples}
\label{sec:numerical_examples}
\subsection{Problem instance of Section~\ref{subsec:counter_ex} -- revisited}
\label{subsection:ex1}
\begin{figure}[!htb]
	\centering
	\includegraphics[width=0.58\linewidth]{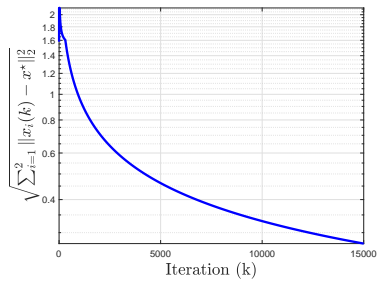}
	\caption{Evolution of $\sqrt{\sum_{i = 1}^2 \| x_i(k) - x^\star \|_2^2}$ for~\eqref{eq:problem_matrices}, where $ (x_i(k))_{k \in \mathbb{N}}, i = 1, 2 $, are the iterates generated by Algorithm~\ref{algo:distributed_algorithm}.}
		\label{fig:Ex1}
\end{figure}

 We revisit the two-agent problem in~\eqref{eq:problem_matrices}, for which the iterative scheme in~\eqref{eq:mod_Duchi} is not guaranteed to converge, and apply this time our algorithm. Note that the optimal solution of~\eqref{eq:problem_matrices} is given by
 \[ 
 x^\star = \mathcal{P}_{[0.5,1]^2} \left[  -\frac{1}{8} Q^{-1} (q_1 + q_2) \right] = \begin{bmatrix}
 0.5\\1
 \end{bmatrix}
 \]
 where $ \mathcal{P}_{[0.5,1]^2} [\cdot] $ represents the projection onto the feasible set of problem~\eqref{eq:problem_matrices}. Pictorially $ x^\star $ is shown in Figure~\ref{fig:Counter_ex}. To illustrate the convergence properties of Algorithm~\ref{algo:distributed_algorithm} we monitor the evolution of $
 \sqrt{\sum_{i = 1}^2 \| x_i(k) - x^\star \|^2_2}$, 
 where $ (x_i(k))_{k \in \mathbb{N}} $, $ i = 1, 2 $, are the iterates generated by Algorithm~\ref{algo:distributed_algorithm}. We use $ c(k) = \frac{1}{\sqrt{k+1}} $ similarly to~\cite{DAW:12}, $ A = \frac{1}{2} \mathbf{1} \mathbf{1}^T $ and $ x_i(0) = \hat{x}_i^\star $, where $ \hat{x}^\star_i $, $ i = 1, 2 $, are defined in Section~\ref{subsec:counter_ex}. Observe that our initial condition is the same as in Proposition~\ref{prop:Counter_ex}. In contrast, as shown in Figure~\ref{fig:Ex1}, the iterates generated by Algorithm~\ref{algo:distributed_algorithm} converge to the optimal solution of~\eqref{eq:problem_matrices}.

\subsection{Example~2: robust linear regression}
\label{subsection:ex2}

We consider the problem of estimating an unknown (but deterministic) vector $ x \in \mathbb{R}^n $ from $ m $ noisy measurements $ y_i $ by means of the linear model
\[ 
y_i = b_i^T x  + v_i, \quad i = 1, \ldots, m,
\]
with $ b_i \in \mathbb{R}^n $, and $ v_i $ are independent random variables drawn from a Laplacian distribution, that is, for each $ i $ the density of $ v_i $ is given by $ h_{v_i}(z) = \frac{1}{2a} \exp^{-|z|/a} $, for all $ z \in \mathbb{R} $.
A common strategy is to impose a norm constraint of the form $  \| x \|_2 \leq c $, for some $ c > 0 $, to reflect some prior knowledge on the unknown vector $ x $, and solve the second order conic program 
\begin{equation}
\hat{x} \in \argmin_{\| x \|_2 \leq c} \|y - B x \|_1.
\label{eq:l1_regression}
\end{equation}
Typically,~\eqref{eq:l1_regression} is referred to as robust regression in the literature, as the $ \ell_1 $-norm penalises relatively less outliers than other convex metrics (e.g., quadratic penalties). In our set-up, we consider the case where data are collected locally and agents are not willing to share their measurements with a central processing unit. 

Observe that~\eqref{eq:l1_regression} has the format of~\eqref{eq:main_pro} by setting  $ X_i = X = \{ x \in \mathbb{R}^n: \| x \|_2 \leq 5 \} $ and $ f_i(x) = | y_i - b_i^T x | $, $ i = 1, \ldots, m $. Moreover, the constraint sets $ X_i $ and the objective functions $ f_i $, $ i = 1, \ldots, m $, trivially satisfy Assumption~\ref{assump:Convexity_Compt}. Hence, we can apply the proposed scheme to obtain a solution to~\eqref{eq:l1_regression}. We consider $ m= 30 $ and $ n = 4 $ and generate $ y $ independently from a standard Gaussian distribution, and matrix $ B $ from a uniform distribution with support $ [0,1] $.

We solve \eqref{eq:l1_regression} in a distributed manner, and compare Algorithm \ref{algo:distributed_algorithm} with the one proposed in \cite{DAW:12} under four different network connectivity structures: $ i) $ complete network graph (which corresponds to the centralised version of the problem); $ ii) $ line network graph; $ iii) $ sparse network graph with sparsity degree $ d = 0.3 $; $ iv) $ sparse network graph with sparsity degree $ d = 0.8 $. We say that a network with $ m $ agents has a sparsity degree $ d \in (0,1) $ if the number of connections among the network nodes is given by $ d m^2 $, where $ m^2 $ indicates the number of connections of a complete graph. 
 
\begin{figure}[!htb]
	\centering
	\begin{psfrags}
		\psfrag{distav1T}[1][1][1][1]{$ $}
		\psfrag{xx}[1][1][1][0]{$ \mathrm{Iterations} $}
		\psfrag{yy}[1][1][1][0]{\hspace{-0.5cm}$  $}
		\includegraphics[width=0.68\linewidth]{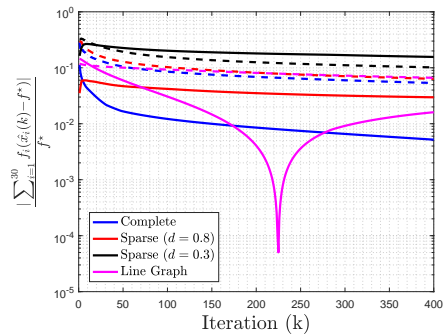}
	\end{psfrags}
	\caption{Evolution of $ \frac{|\sum_{i = 1}^{30} f_i(x_i(k)) - f^\star|}{f^\star} $ for Algorithm~\ref{algo:distributed_algorithm} (solid lines) and the one in~\cite{DAW:12} (dashed lines) when applied to the robust regression problem given by~\eqref{eq:l1_regression}. The different colours correspond the different network connectivities.}
		\label{fig:Ex2}
\end{figure}

 We assess the performance of Algorithm~\ref{algo:distributed_algorithm} for each of the aforementioned networks in Figure~\ref{fig:Ex2}. Solid lines correspond to Algorithm~\ref{algo:distributed_algorithm}, whereas dashed lines correspond to the algorithm proposed in~\cite{DAW:12}. Different colours correspond to the different network connectivities. For each case, we monitor the evolution of $ \frac{|\sum_{i = 1}^{30} f_i(x_i(k)) - f^\star|}{f^\star} $, where $ f^\star $ is the optimal value of~\eqref{eq:l1_regression}. The proposed scheme exhibits similar and often favourable performance with the one in~\cite{DAW:12}, in particular for cases where the underlying graph is not sparse. It should be noted, however, that Algorithm~\ref{algo:distributed_algorithm} possesses more general convergence properties, i.e., the proposed scheme is guaranteed to converge under non-identical local sets.
 
Note that due to the fact that Algorithm~\ref{algo:distributed_algorithm} requires two rounds of communication per iteration, the results presented in Figure \ref{fig:Ex2} should be rescaled by a factor of two if we use communication rounds instead of the iteration index.
  

\subsection{Example~3: $ \ell_2 $ linear regression with regularisation}
\label{subsec:Ex3}

In this example, we consider a variation of the regression problem where we assume $ v_i $, $ i = 1, \ldots, m $, to be independent and Gaussian, i.e., the density function is given by $ h_{v_i}(z) = \frac{1}{\sqrt{2\pi}} e^{-\frac{z^2}{2}}  $, for all $ z \in \mathbb{R} $, for all $ i = 1, \ldots, m $, and we assume that $ x $ is sparse. A common relaxation of this problem is to choose the maximum likelihood estimator $ \hat{x} $ such that 
\begin{equation}
\begin{aligned}
\hat{x} = \argmin_{x \in X} \| y - B x \|_2^2 + \lambda \|x \|_1,
\end{aligned}
\label{eq:reg_l2_l1_norm}
\end{equation}
where $ X $ can be interpreted as a set including prior beliefs, e.g., $ \|x \|_2 \leq c $ or $ \underline{x} \leq x \leq \bar{x} $ for some vectors $ \underline{x}, \bar{x} \in \mathbb{R}^n $. The estimator $ \hat{x} $ obtained by solving~\eqref{eq:reg_l2_l1_norm} depends on the value of the parameter $ \lambda$. In fact, the larger the value of $\lambda$, the worse the performance is in terms of the error and the sparser the obtained solution is.  

In this example, we aim to verify the performance of Algorithm~\ref{algo:distributed_algorithm} under step size choices $ c(k) \propto \frac{1}{k+1} $ and a time-varying communication network. Similar to the previous example, the vector $ y $ is generated according to a standard normal distribution and matrix $ B $ from a uniform distribution on the interval $ [0,1] $. We assume $ m > n $ and consider the case where agents possess private, local information, encoded by $ X_i = [\underline{x}_i,\bar{x}_i] $ $ i = 1, \ldots, m $, such that $ X = \cap_{i = 1}^m X_i = [\underline{x},\bar{x}] $. 
%

The algorithm presented in~\cite{DAW:12} does not necessarily converge in the set-up of problem~\eqref{eq:reg_l2_l1_norm}, as we have different constraint sets per agent.  We thus compare our algorithm against the one proposed in~\cite{MFGP:18}, which converges under similar conditions but does not leverage on subgradient averaging. This allows us to assess the impact of averaging subgradients on practical convergence.

\begin{figure}[!htb]
	\centering
	\includegraphics[width=0.68\linewidth]{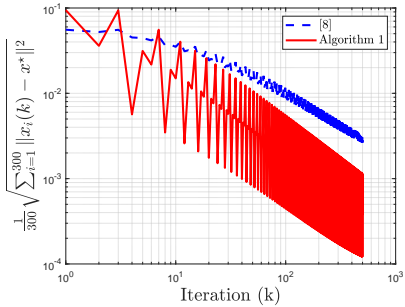}
	\caption{Evolution of the average distance to the optimal solution given by $ \frac{\sqrt{\sum_{i = 1}^{300} \| x_i(k) - x^\star \|^2_2}}{300} $ for Algorithm~\ref{algo:distributed_algorithm} (solid-red line) and that of~\cite{MFGP:18} (dashed-blue line).}
	\label{fig:Ex3b}
\end{figure}

%
%
%
We now investigate the behaviour of the proposed algorithm in the presence of time-varying communication networks. To this end, we set $ m = 300$ and $ n = 10 $, and generate four network configurations with different sparsity patterns, alternating cyclically among these. We also set $ c(k) = \frac{0.2}{k+1} $ for both Algorithm~\ref{algo:distributed_algorithm} and the one in~\cite{MFGP:18}. Figure~\ref{fig:Ex3b} shows the evolution for the average distance to the optimal solution for Algorithm~\ref{algo:distributed_algorithm} (solid-red line) and the one in~\cite{MFGP:18} (dashed-blue line). We observe that Algorithm~\ref{algo:distributed_algorithm} consistently outperforms the one proposed in~\cite{MFGP:18}; this is mainly due to the sub-gradient averaging step of Algorithm~\ref{algo:distributed_algorithm}. 

\section{Conclusion}
\label{sec:conclusion}

In this paper we proposed a subgradient averaging algorithm for multi-agent optimisation problems involving non-differentiable objective functions and different constraint sets per agent. For this set-up we showed by means of a geometric construction that available schemes involving subgradient averaging cannot be used. For the proposed scheme we showed convergence of the algorithm iterates to some minimiser of a centralised problem counterpart. Moreover, we have also established a convergence rate under a particular choice for the underlying step size.
The performance of our approach was illustrated by means of several numerical examples, quantifying also the  improvement in terms of practical convergence with respect to other algorithms that are not based on (sub)gradient exchange.

Future work will concentrate towards replacing the diminishing step size employed by our approach with a constant one, showing convergence rates to a neighbourhood of the set of optimal solutions. A more detailed study on the communication requirements, and an investigation on how we could reduce the two rounds of communication required by the proposed algorithm is also a topic of current work.

\section{Appendix}
\label{sec:appendix}

\subsection{Proof of Lemma~\ref{lemma:bounded_subgradient}}
\label{app:Prof_lemma_subgrad}

We start by proving item $ i) $. Consider the continuous mapping $ \phi: \mathbb{R}^m \times \prod_{i = 1}^{m} \mathbb{R}^n \rightarrow \mathbb{R}^n$, defined as $ 
\phi(\gamma,x_1,\ldots,x_m) = \sum_{i = 1}^m \gamma_i x_i,$ where $ \gamma = (\gamma_1,\ldots,\gamma_m) $ belongs to the simplex in $ \mathbb{R}^m $, denoted by $ \Gamma $. Consider $ K = \phi(\Gamma,\prod_{i = 1}^m X_i) $, and note that $ K $ is compact, as it is the image of the compact set $ \Gamma \times \prod_{i = 1}^m X_i $ under the continuous map $ \phi $. Moreover, note that by definition we have $ K \subseteq \mathrm{conv}(\cup_{i = 1}^m X_i), $ as any element in $ K $ is a convex combination of elements in $ \cup_{i = 1}^m X_i $. To conclude the argument, we need to show that $ \mathrm{conv}(\cup_{i = 1}^m X_i) \subseteq K $. To this end, it suffices to show that $ K $ is a convex set, due to the fact that the convex hull is the smallest convex set containing a given set. Let $ z, w \in K $, i.e.,  $ z = \sum_{i = 1}^m \gamma_i z_i $ and $  w = \sum_{i = 1}^m \beta_i w_i $, with $ z_i, w_i \in X_i $, and $ \gamma = (\gamma_1,\ldots,\gamma_m) , \beta = (\beta_1,\ldots,\beta_m) \in \Gamma$. Fix an $ \alpha \in (0,1) $, and note that
$
		 \alpha z + (1-\alpha) w = \sum_{i = 1}^m (\alpha \gamma_i + (1-\alpha)\beta_i) x_i,
		 $ where $ x_i = c_i z_i + (1-c_i) w_i \in A_i, $ with $ c_i = \frac{\alpha \gamma_i}{\alpha \gamma_i + (1-\alpha) \beta_i} $. 
		 
		 Since $ x_i \in A_i $ due to convexity of $ A_i $ and $ \alpha \gamma + (1-\alpha) \beta \in \Gamma$, we conclude that $ \alpha z + (1-\alpha)w \in K $ for any $ \alpha \in (0,1) $, thus showing that $ K $ is a convex set. This implies then that $ K = \mathrm{conv}(\cup_{i = 1}^m X_i) $ as we have established that $ K \subseteq \mathrm{conv}(\cup_{i = 1}^m X_i) $ and $ \mathrm{conv}(\cup_{i = 1}^m X_i) \subseteq K $. Since $ K $ was shown to be compact, we have that $ \mathrm{conv}(\cup_{i = 1}^m X_i) $ is also compact. This concludes the proof of item $ i). $ An alternative proof can be found at~\cite[Prop. 1.2.2]{Ber:09}. The proof of item $ ii) $ follows from Proposition~5.4.2, p. 186, in~\cite{Ber:09}, and is omitted for brevity. This concludes the proof of the lemma.

\subsection{Sufficient condition for Assumption~\ref{assump:Convexity_Compt}, item $ iii) $.}
\label{app:suff_bounded_subgrad}

The goal of this subsection is to provide a sufficient condition for Assumption~\ref{assump:Convexity_Compt}, item $ iii) $. The subsequent arguments can be found in standard optimisation books, such as~\cite[Theorem 24.7]{Roc:72}; however we present here a more direct proof.
\begin{assump}Let $ X_i $, $ i = 1, \ldots, m $, be the level sets of problem~\eqref{eq:main_pro} and $ \mathrm{dom} f $ the domain of $ f $. We suppose that:
	
	\begin{itemize} 
		\item[$ i) $] The distance between the set $ \cup_{i = 1}^m X_i $ and the complement of the interior of the domain of $ f $ (which is closed and convex) is strictly greater than zero, i.e.,
		\begin{equation*}
		\begin{aligned}
		\mathrm{dist}(&\cup_{i = 1}^m X_i, \big(\mathrm{int} (\mathrm{dom} f)\big)^c)  \\ &= \inf_{\substack{x \in \cup_{i = 1}^m X_i,\\ y \in \big(\mathrm{int} (\mathrm{dom} f)\big)^c }} \| x - y \|_2^2> 0.
		\end{aligned}
		\end{equation*}
		\item[$ ii) $] $ X_i \subset \cap_{i = 1}^m \mathrm{int}(\mathrm{dom} f_i) $ for each $ i = 1, \ldots, m $.
	\end{itemize}
	\label{assump:bounded_subgrad}
\end{assump}

As a consequence of Assumption~\ref{assump:bounded_subgrad}, and since $  \mathrm{dom} f = \cap_{i = 1}^m \mathrm{dom} f_i $, $ \mathrm{ri}(\mathrm{dom} f) = \cap_{i = 1}^m \mathrm{ri}(\mathrm{dom} f_i) $ and $ \mathrm{ri}(\mathrm{dom} f_i) \subset \mathrm{dom} f_i $ we have that the subdifferential $ \partial f(x) $  is non empty for each $ x \in \cap_{i = 1}^m X_i $, as by item $ ii) $ of Assumption~\ref{assump:bounded_subgrad} every feasible solution of~\eqref{eq:main_pro} belongs to the interior of the domain of $ f $. Furthermore, $ \partial f(x) $ is compact by~\cite[Proposition~5.4.1]{Ber:09} since the affine hull of $ \mathrm{dom} f $ has dimension $ n $ due to Assumption~\ref{assump:Convexity_Compt}, item $ ii) $. 
We use this fact to show that $\cup_{x \in \mathrm{conv}(\cup X_i)} \partial f(x) $ is a bounded set, that is, $ \| g \|_2 \leq L $, where $ g \in \partial f(x) $ for any $ x \in \cup_{i = 1}^m X_i $. This result is formally stated in the next lemma.
\begin{lemma}
	Under Assumptions~\ref{assump:Convexity_Compt}, items $ i) $ and $ ii) $ , and~\ref{assump:bounded_subgrad}, we have that the set $ \cup_{x \in \mathrm{conv}(\cup X_i)} \partial f(x) $ is non-empty and bounded.
\end{lemma}

\begin{pf}
	
	The proof of the lemma relies on Assumption~\ref{assump:bounded_subgrad}, item $ ii) $, that is, $ X_i \subset \cap_{j = 1}^m \mathrm{ri}(\mathrm{dom} f_j) $, for all $ i = 1, \ldots, m $. This implies that $ \mathrm{conv}(\cup_{i = 1}^m X_i) \subset \cap_{j = 1}^m \mathrm{ri} (\mathrm{dom} f_j) $, as $ \cap_{j = 1}^m \mathrm{ri} (\mathrm{dom} f_j) $ is convex and contains $ \cup_{i = 1}^m X_i $. Suppose, by contradiction, that $ \cup_{x \in \mathrm{conv}(\cup X_i)} \partial f(x) $ is unbounded. Then there exists a sequence $ (x_k)_{k \in \mathbb{N}} \subset \mathrm{conv} (\cup_{i = 1}^m X_i) $ such that $ (g_k)_{k \in \mathbb{N}} $, with $ g_k \in \partial f(x_k) $, satisfies  $ 
	\| g_k \|_2 < \| g_{k + 1} \|_2, \quad \forall~k~\in~\mathbb{N}.
	$
	
	Notice that $ x_k \in \cap_{i = 1}^m \mathrm{int}(\mathrm{dom} f_i) $ by Assumption~\ref{assump:bounded_subgrad}, item $ ii) $.  By item $ i) $ of Assumption~\ref{assump:bounded_subgrad}, we can construct a sequence $ (\beta_k)_{k \in \mathbb{N}} $ such that $
	x_k + \beta_k d_k \in  \cap_{i = 1}^m \mathrm{dom} f_i.$
	with $ d_k = g_k/\| g_k \|_2 $. Let $ \beta = \inf_{k \in \mathbb{N}} \beta_k $ and notice that $ \beta > 0 $ (i.e., it is bounded away from zero) due to Assumption~\ref{assump:bounded_subgrad}, item $ i).$ By the definition of $ g_k $ we have that 
	\begin{equation}
	\frac{f(x_k + \beta d_k) - f(x_k)}{\beta} \geq \| g_k \|_2, \quad \forall~k~\in~\mathbb{N}.
	\label{eq:lemma_1}
	\end{equation}
	As inequality~\eqref{eq:lemma_1} is valid for all $ k \in \mathbb{N} $, we take the limit superior on both sides to obtain
	\small
	\begin{equation}
	\limsup_{k \rightarrow \infty} \| g_k \|_2 \leq \limsup_{k \rightarrow \infty} \frac{f(x_i + \gamma d_k) - f(x_k)}{\gamma} < \infty,
	\label{eq:lemma_2}
	\end{equation}
	\normalsize
	where the right-hand side of~\eqref{eq:lemma_2} is finite as the sequences $ (x_k)_{k \in \mathbb{N}} $ and $ (d_k)_{k \in \mathbb{N}} $ are bounded (notice that $ d_k $ is a normalised subgradient), and since $ f $ is continuous on its domain ($ f $ is convex). This establishes a contradiction, as we assumed $ (g_k)_{k \in \mathbb{N}} $ were unbounded, thus concluding the proof of item $ ii) $. 
	
\end{pf}

\subsection{Proof of Proposition~\ref{prop:Counter_ex}}
\label{app:Proof_prop_counter}
The proof is based on an induction argument.
\subsubsection*{Base case} We show that $ z_i(1)^T(\xi - \hat{x}_j^\star) \geq 0, $ for all $ \xi \in X_j $, for all $ i, j = 1, 2 $, and also that $ x_i(1) = \hat{x}^\star $, for all $ i = 1, 2 $. Consider the inequalities 
	\begin{align}
	&\nabla f_1(\hat{x}_1^\star)^T (\xi - \hat{x}_i^\star) \geq 0, \nonumber \\&\nabla f_2 (\hat{x}_2^\star)^T  (\xi - \hat{x}_i^\star) \geq 0, \quad \forall \xi \in X_i, \quad i = 1, 2.
	\label{eq:first_ineq_proof}
	\end{align}
	Fix $ i = 1  $. The first inequality in~\eqref{eq:first_ineq_proof} holds due to optimality of~$ \hat{x}^\star_1 $~\cite{Ber:09}. To show the second inequality observe that $ 
	\nabla f_2 (\hat{x}_2^\star) = \left[ 
	13.68,
	-3.94
	\right]^T,$
	and that $ \xi - \hat{x}_1^\star = [a_1,a_2]^T $ with $ a_1 \geq 0 $ and $ a_2 \leq 0 $, for all $ \xi \in X_1 $. 
	
	 Since $ \nabla f_1 (\hat{x}_1^\star) =\left[ 
      12, -4 \right]^T,$
	using a symmetric argument we show that 
	\begin{align}
	&\nabla f_2(\hat{x}_2^\star)^T (\xi - \hat{x}^\star_2) \geq 0, \nonumber \\ &\nabla f_1(\hat{x}_1^\star)^T (\xi - \hat{x}^\star_2) \geq 0, \quad \forall \xi \in X_2.
	\label{eq:second_ineq}
	\end{align}
	By~\eqref{eq:aveg_step_mod_Duchi}, and under our choice for $ A $,
	\begin{equation}
	   z_i(1) = \frac{1}{2} \Big(  \nabla f_1(\hat{x}_1^\star) + \nabla f_2(\hat{x}_2^\star) \Big) + \nabla f_i(\hat{x}^\star_i),   
	 \label{eq:z_def_proof}
	\end{equation}
	for $ i = 1, 2$, hence inequalities~\eqref{eq:first_ineq_proof} and~\eqref{eq:second_ineq} imply that $ z_i(1)^T (\xi - \hat{x}^\star_j) \geq 0,~\forall \xi \in X_j $, for all $ i, j = 1,2 $. 
	
	We will now prove that $ x_i(1) = \hat{x}^\star_i $, for $ i = 1, 2 $. Fix $ i = 1 $. Since $ z_1(1)^T \xi  + \frac{2}{c(k)} \| \xi \|_2^2$ is strictly convex, there is a unique point satisfying
	\begin{equation}
	\Big( z_1(1) + 2 x_1(1)\Big)^T (\xi - x_1(1)) \geq 0, \quad \forall \xi \in X_1,
	\label{eq:forth_ineq_pro}
	\end{equation}
	where $ (z_1(1) + 2 x_1(1)) $ is the gradient of the objective function in~\eqref{eq:x_update_mod_Duchi} evaluated at $ x_1(1) $, with $ c(1) = 1 $.  Therefore, it suffices to show that 
		\begin{equation}
	\Big( z_1(1) + 2 \hat{x}^\star_1\Big)^T (\xi - \hat{x}^\star_1) \geq 0, \quad \forall \xi \in X_1.
	\label{eq:fifth_ineq_prop}
	\end{equation}
	By substituting~\eqref{eq:problem_matrices} into~\eqref{eq:z_def_proof}, we observe that $ 
	z_1(1) + 2\hat{x}^\star_1 =\left[ 
	22.8414,
	-5.9708
	\right]^T,
	$
	and due to the structure of $ \xi - \hat{x}^\star_1 $,~\eqref{eq:fifth_ineq_prop} holds, thus proving that $ x_1(1) = \hat{x}^\star_1 $. A symmetric argument yields that $ x_2(1) = x^\star_2 $.
	
	\subsubsection*{Induction hypothesis} Assume that $ z_i(k)^T (\xi-\hat{x}_j^\star) \geq 0$ for all $ \xi \in X_j $, for $ i,j = 1,2 $, and that $ x_i(k) = x^\star_i $ for $ i = 1,2 $. We aim to show that the aforementioned relations remain true for the step $ k+1. $
	
	\subsubsection*{Proof for iteration $ k + 1 $} Fix $ i = 1 $. Following a similar reasoning with the base case, observe that  $ x_1(k+1) = x^\star_1  $ if
	\begin{equation}
	\left[ z_1(k+1) + \frac{2}{c(k)} \hat{x}^\star_1 \right]^T (\xi - \hat{x}^\star_1) \geq 0, \quad \forall \xi \in X_1.
	\label{eq:third_eq}
	\end{equation}
	 As the sequence $ (z_i(k))_{k \in \mathbb{N}} $ is generated by~\eqref{eq:aveg_step_mod_Duchi}, we propagate the dynamical system in~\eqref{eq:aveg_step_mod_Duchi} by $ k+1 $ steps to obtain
	\[ 
	z_i(k+1) = \frac{1}{2} \Big( \nabla f_1 (\hat{x}_1^\star) 
	+  \nabla f_2 (\hat{x}_2^\star)  \Big)(k+1) + \nabla f_1 (\hat{x}^\star_1),
	\]
	where we have used the fact that $ A = \frac{1}{m} \mathbf{1} \mathbf{1}^T $ and $ c(k) = \frac{1}{\sqrt{k+1}} $. A sufficient condition for equation~\eqref{eq:third_eq} to hold is that 
	\begin{align}
	\Bigg[ \frac{1}{2} \Big( \nabla &f_1 (\hat{x}_1^\star) 
	+  \nabla f_2 (\hat{x}_2^\star)  \Big) (k+1) \nonumber  \\   &+ 2 \hat{x}^\star_1 \sqrt{k+1} \Bigg]^T (\xi - \hat{x}^\star_1) \geq 0, 
	\quad \forall \xi \in X_1,
	\label{eq:forth_ineq}
	\end{align}
	since $ \nabla f_1(\hat{x}^\star_1)^T (\xi - \hat{x}^\star_1)\geq 0 $ by optimality of $ \hat{x}^\star_1 $.
	Recall that $ (\xi - \hat{x}_1^\star) = [a_1,a_2] $ with $ a_1 \geq 0 $ and $ a_2 \leq 0 $ for all $ \xi \in X_1 $. To prove~\eqref{eq:forth_ineq} we will show that the left-most vector in the same equation can be written as $ [b_1,b_2] $ for some $ b_1 \geq 0 $ and $ b_2 \leq 0 $. To achieve this, notice that $ k + 1 \geq \sqrt{2} \sqrt{k+1}  $, for all $ k \geq 1 $, and let $ e_i $ denote the unit vector with $ 1 $ in the $ i$-th position, $  i = 1, 2 $. We then have that	
	\begin{align}
	e_1^T &\left[ \frac{1}{2} \Big( \nabla f_1 (\hat{x}_1^\star) 
	+  \nabla f_2 (\hat{x}_2^\star)  \Big) \right] (k+1) \nonumber \nonumber \\
	&\geq e_1^T \left[ \frac{\sqrt{2}}{2} \Big( \nabla f_1 (\hat{x}_1^\star) 
	+  \nabla f_2 (\hat{x}_2^\star)  \Big) \right] \sqrt{k+1},
	\label{eq:e_1_lb_proof_prop}
	\end{align}
	and
	\begin{equation}
	2e_2^T   \hat{x}^\star_1 \sqrt{k+1} \leq \sqrt{2} e_2^T \hat{x}_1^\star (k+1),
	\label{eq:e_2_ub_proof_prop}
	\end{equation}
	since the first component of the averaged gradient and the second component of $ \hat{x}^\star_1 $ are both positive. Therefore, for all $ k \in \mathbb{N} $,
	\small
	\begin{align}
	b_1 \geq 16.1604 \sqrt{k+1} > 0, \quad
	b_2 \leq -2.5566 (k+1) < 0.
	\label{eq:seventh_eq}
	\end{align}
	\normalsize
	Inequalities~\eqref{eq:e_1_lb_proof_prop},~\eqref{eq:e_2_ub_proof_prop} and~\eqref{eq:seventh_eq}, together with the structure of $ \xi - \hat{x}_1^*  $, imply that~\eqref{eq:forth_ineq} holds, so we can conclude that $ x_1(k+1) = \hat{x}^\star_1 $. A symmetric argument shows that $ x_2(k+1) = \hat{x}^\star_2 $. 
	
	To complete the proof it remains to show that $ z_i(k+1)^T (\xi-\hat{x}_j^\star) \geq 0$ for all $ \xi \in X_j $, for all $ i,j = 1,2 $, where 
		$
	z_i(k+1) = \frac{1}{2} \Big( z_1(k) + z_2(k) \Big) + \nabla f_i(x_i(k)),   
	$
	due to~\eqref{eq:aveg_step_mod_Duchi} and our choice for $ A $. By our induction hypothesis, $ z_i(k) (\xi - \hat{x}^\star_j) \geq 0$, for all $ i,j = 1,2 $, hence it suffices to show that $  \nabla f_i(x_i(k))^T (\xi - \hat{x}^\star_j) \geq 0,~\forall \xi \in X_j,~\forall  i = 1, 2. $
	Since $ x_i(k) = \hat{x}^\star_i $ for $ i = 1,2 $, due to our induction hypothesis, the claim follows from~\eqref{eq:first_ineq_proof} and~\eqref{eq:second_ineq}, thus concluding the proof.

\subsection{Auxiliary Lemmas for the proofs of Theorem~\ref{theo:convergence_algo_1} and~\ref{theo:convergence_rate}.}
\label{app:Auxiliary_results}
Let 
\begin{equation}
v(k) = \frac{1}{m} \sum_{i = 1}^m x_i(k),
\label{eq:v_def}
\end{equation}
be the average of the agents' estimates at time $ k $. Since this quantity might not necessarily belong to the feasible set $ \cap_{i = 1}^m X_i $, we define 
\begin{equation}
\bar{v}(k) = \frac{\rho}{\epsilon(k) + \rho} v(k) +  \frac{\epsilon(k)}{\epsilon(k) + \rho} \bar{x}, 
\label{eq:def_v_bar}
\end{equation}
where $ \bar{x} $ is a point in the interior of the feasible set (which is non-empty by Assumption~\ref{assump:Convexity_Compt}, item $ ii) $), $ \rho > 0 $ is such that the $ 2 $-norm ball of centre $ \bar{x} $ and radius $ \rho $ is contained in $ \cap_{i = 1}^m X_i $, and $ \epsilon(k) = \sum_{i = 1}^m \mathrm{dist}(v(k),X_i) $. As shown in~\cite{NOP:10}, $ \bar{v}(k) \in  \cap_{i = 1}^m X_i $, for all $ k \in \mathbb{N} $. We also define $ e_i(k+1) = x_i(k+1) - z_i(k) $, and note that the $ z_i $-update in Algorithm~\ref{algo:distributed_algorithm} can be written as
\begin{equation}
x_i(k+1) = \sum_{j = 1}^m [A(k)]_j^i x_j(k) + e_i(k+1).
\label{eq:error_def}
\end{equation}

\begin{lemma}
	The following relations hold.
	\begin{itemize}
		\item[$ i)  $] Let $ (x_i(k))_{k \in \mathbb{N}} $, $ i = 1, \ldots, m $, be the sequences generated by Algorithm~\ref{algo:distributed_algorithm}, and $ (v(k))_{k \in \mathbb{N}} $ and $ (\bar{v}(k))_{k \in \mathbb{N}} $ defined by~\eqref{eq:v_def} and~\eqref{eq:def_v_bar}, respectively. Under Assumption~\ref{assump:Convexity_Compt}, we have that for all $ k \geq 0 $, 
		\small
		\[
		\sum_{i = 1}^m \| x_i(k+1)- \bar{v}(k) \|_2 \leq \mu \sum_{i = 1}^m \| x_i(k) - v(k) \|_2, 
		\]
		\normalsize
		where $ \mu = \frac{2}{\rho} m D + 1 $, and $ D $ is the diameter of the set $ \cup_{i = 1}^m X_i $ (which is well-defined by Lemma~\ref{lemma:bounded_subgradient}, item $ i) $).
		\item[$ ii) $] Let $ (x_i(k))_{k \in \mathbb{N}} $, $ i = 1, \ldots, m $, and $ (v(k))_{k \in \mathbb{N}} $ be as in item $ i) $. Under Assumption~\ref{assump:Strongly_graph}, we have that for all $ i = 1, \ldots, m $, for all $ k \geq 0  $,
		\small
		\begin{equation*}
		\begin{aligned}
		\| x_i&(k+1) - v(k+1) \|_2\leq \lambda q^{k} \sum_{j = 1}^m \| x_j(0) \|_2 \\& + \| e_i(k+1) \|_2 + \sum_{r = 0}^{k - 1} \lambda q^{k-r-1} \sum_{j = 1}^m \| e_j (r+1) \|_2 \\ &+  \frac{1}{m} \sum_{j = 1}^m \| e_j(k+1) \|_2,
		\end{aligned}
		\end{equation*}
		\normalsize
		where $ \lambda = 2(1 + \eta^{-(m-1)T})/(1-\eta^{(m-1)T}) \in \mathbb{R}_+$ and $ q = (1-\eta^{(m-1)T})^{\frac{1}{(m-1)T}} \in (0,1)$.
		\item[$ iii) $] Given a non-increasing and non-negative sequence $ (c(k))_{k \in \mathbb{N}} $, and a scalar $ \bar{L} > 0 $, we have that
		\small
		\begin{equation*}
		\begin{aligned}
		2 \bar{L} &\sum_{k = 0}^{N} c(k) \sum_{i = 1}^{m} \| x_i(k+1) - \bar{v}(k + 1) \|_2 \\ &< \beta_1 \sum_{k = 0}^{N}\sum_{i = 1}^{m} \| e_i(k+1) \|_2^2 + \beta_2 \sum_{k = 0}^{N} c(k)^2 + \beta_3, 
		\end{aligned}
		\end{equation*}
		\normalsize
		where $ \beta_1 \in (0,1) $, and $ \beta_2 $ and $ \beta_3 $ are positive constants. 
	\end{itemize}
	\label{lemma:literature_results}
\end{lemma}
\begin{pf}
	The proof of item $ i) $ is presented in~\cite[Lemma~1]{MFGP:18}. For item $ ii) $, see~\cite[Lemma~2]{MFGP:18}. Finally, the proof of item $ iii) $ follows the line of~\cite[Lemma~3]{MFGP:18}. 
\end{pf}

Observe that the values of $ \lambda $ and $ q $ in Lemma~\ref{lemma:literature_results}, item $ ii) $, depend on the parameter $ T $ that characterises the uniform bound in Assumption~\ref{assump:Strongly_graph}, item $ i) $; and on $ \eta $, the lower bound for the elements of $ A(k) $, Assumption~\ref{assump:Strongly_graph}, item $ ii) $. The following lemma is instrumental for the proof of Theorem~\ref{theo:convergence_rate}. In particular, Lemma~\ref{lemma:main_ineq}, item $ ii) $, constitutes a non-trivial extension of the result in \cite{MFGP:18}, allowing some sequences to be iteration-varying.

\begin{lemma}
	Let $ (x_i(k))_{k \in \mathbb{N}}, (z_i(k))_{k \in \mathbb{N}} $ and $ (d_i(k))_{k \in \mathbb{N}} $, $ i = 1, \ldots, m $, be the sequences generated by Algorithm~\ref{algo:distributed_algorithm}, and $ x^\star $ by any optimal solution of~\eqref{eq:main_pro}. Under Assumptions~\ref{assump:Convexity_Compt} and~\ref{assump:Strongly_graph}, we have that:
	\begin{itemize}
		\item[$ i) $]  For all $ k \in \mathbb{N} $, 
		\small
		\begin{align}
		&2c(k) \sum_{i = 1}^m d_i(k)^T (x_i(k+1) - x^\star) + \sum_{i = 1}^m \| e_i(k+1) \|_2^2 \nonumber\\
		&+ \sum_{i = 1}^m \| x_i(k+1) - x^\star \|_2^2   \leq \sum_{i = 1}^m \| x_i(k)- x^\star \|_2^2.
		\label{eq:main_ineq}
		\end{align}
		\normalsize
		\item[$ ii) $] For any $ \beta_1 \in (0,1) $, there exist sequences $ (\alpha_1(k))_{k \in \mathbb{N}} $ and $ (\alpha_2(k))_{k \in \mathbb{N}} $ such that, for all $ k \in \mathbb{N} $, $ \alpha_1(k) \in (0,1) $, $ \alpha_2(k) \in (0,1) $, $ 1-\beta_1-\alpha_1(k) -\alpha_2(k) \geq 0 $ and 
		\small
		\begin{align}
		&2 \sum_{k = 0}^N c(k) \sum_{i = 1}^m (f_i(\bar{v}(k+1)) - f_i(x^\star)) \nonumber \\  &+  \sum_{k = 0}^N (1-\alpha_1(k) - \alpha_2(k) - \beta_1) \sum_{i = 1}^m  \| e_i(k+1) \|_2^2  \nonumber   \\ & + \sum_{k = 0}^N \sum_{i = 1}^m \| x_i(k+1) - x^\star \|_2^2    \leq \sum_{k = 0}^N \sum_{i = 1}^m \| x_i(k)- x^\star \|_2^2 \nonumber \\ &+  \sum_{k = 0}^N \Big( m L^2  \frac{\alpha_1(k) + \alpha_2(k)}{\alpha_1(k) \alpha_2(k)} + \beta_2 \Big) c(k)^2  + \beta_3.
		\label{eq:main_ineq_with_k}
		\end{align}
		\normalsize
	\end{itemize}
	\label{lemma:main_ineq}
\end{lemma}
\begin{pf}
	Item $ i) $: Fix any $ i \in \{1, \ldots, m\}$ and consider the sequence $ (x_i(k))_{k \in \mathbb{N}} $. By optimality of $ x_i(k+1) $ (see Algorithm~\ref{algo:distributed_algorithm}), for any $ \xi \in X_i $,
	\small
		\begin{align}
	d_i(k)^T &x_i(k+1) -\frac{1}{c(k)} (z_i(k) - x_i(k+1))^T x_i(k+1) \nonumber \\
	&\leq d_i(k)^T \xi - \frac{1}{c(k)} (z_i(k) - x_i(k+1))^T \xi,
	\label{eq:proof_ineq_lemma_1}
	\end{align}
	\normalsize
	where $ d_i(k) - \frac{1}{c(k)} (z_i(k) - x_i(k+1)) $ constitutes the gradient of the objective function in the $ x_i-$update of Algorithm~\ref{algo:distributed_algorithm}, evaluated at $ x_i(k+1) $. Fix any optimal solution of~\eqref{eq:main_pro}, $ x^\star \in \cap_{i = 1}^m X_i, $ and consider the following identity
	\small
	\begin{align}
	&-\frac{1}{c(k)} (z_i(k) - x_i(k+1))^T(x_i(k+1) - x^\star)\nonumber \\ &= \frac{1}{2c(k)} \| x_i(k+1) - z_i(k) \|_2^2 + \frac{1}{2 c(k)} \| x_i(k+1) - x^\star \|_2^2 \nonumber \\
	&-\frac{1}{2 c(k)} \| z_i(k) - x^\star \|_2^2.
	\label{eq:proof_2}
	\end{align}
	\normalsize
	Combining~\eqref{eq:proof_2} and~\eqref{eq:proof_ineq_lemma_1} with $ \xi = x^\star $, we obtain
	\small
	\begin{align}
	&d_i(k)^T x_i(k+1) + \frac{1}{2c(k)} \| x_i(k+1) - z_i(k) \|_2^2 \nonumber \\ & \quad \quad \quad+ \frac{1}{2 c(k)} \| x_i(k+1) - x^\star \|_2^2 \nonumber \\ &\leq d_i(k)^T x^\star + \frac{1}{2 c(k)} \| z_i(k) - x^\star \|_2^2 \nonumber \\
	& \leq  d_i(k)^T x^\star + \frac{1}{2 c(k)} \sum_{j = 1}^m [A(k)]_j^i \| x_j(k) - x^\star \|_2^2,
	\label{eq:proof_ineq_lemma_2}
	\end{align}
	\normalsize
	where the last inequality follows from double stochasticity of $ A(k) $ and convexity of $ \| \cdot \|^2 $.
	
	We now multiply both sides of~\eqref{eq:proof_ineq_lemma_2} by $ 2 c(k) $ and sum the result for all $ i = 1, \ldots, m, $ to obtain
	\small
	\begin{align}
	2c(k)&\sum_{i = 1}^m d_i(k)^T x_i(k+1) + \sum_{i = 1}^m \| x_i(k+1) - z_i(k) \|_2^2 \nonumber \end{align} \begin{align} &+ \sum_{i = 1}^m \| x_i(k+1) - x^\star \|_2^2   \leq  2c(k) \sum_{i = 1}^m d_i(k)^T x^\star \nonumber \\
	&+  \sum_{i = 1}^m  \| x_i(k) - x^\star \|_2^2,
	\label{eq:proof_ineq_lemma_3}
	\end{align}
	\normalsize
	where $ \sum_{i = 1}^{m} \sum_{j = 1}^{m} [A(k)]_j^{i} \| x_j(k) - x^\star \|_2^2 = \sum_{i =1}^{m} \| x_i(k) - x^\star\|_2^2 $ by exchanging the order of summation, and due to double stochasticity of $ A(k) $. The result follows from~\eqref{eq:proof_ineq_lemma_3} by recalling that $ e(k+1) = x_i(k+1) - z_i(k) $ and moving the first term in the right-hand side of~\eqref{eq:proof_ineq_lemma_3} to the left one. This concludes the proof of item $ i) $.
	
	Item $ ii) $: Consider the first term in the left-hand side of~\eqref{eq:main_ineq}, and rewrite it as\small
	\begin{align}
	2c(k)& \sum_{i = 1}^m d_i(k)^T (x_i(k+1) - x^\star) = \nonumber  \\ &2c(k) \sum_{i = 1}^m d_i(k)^T (x_i(k+1) - \bar{v}(k+1)) \nonumber \\ + &2c(k) \sum_{i = 1}^m d_i(k)^T (\bar{v}(k+1) - x^\star) 
	\label{eq:ineq_1}
	\end{align}
	\normalsize
	by adding and subtracting $ \bar{v}(k+1) $. We next consider the terms in the right hand-side of~\eqref{eq:ineq_1} separately. First, observe that
	\small
	\begin{align}
	2c(k) &\sum_{i = 1}^m d_i(k)^T (x_i(k+1) - \bar{v}(k+1))\nonumber \\  \geq &-2c(k) L \sum_{i = 1}^m  \| x_i(k+1) - \bar{v}(k + 1) \|_2,
	\label{eq:ineq_2}
	\end{align}
	\normalsize
	by the Cauchy-Schwartz inequality, where $	L =  \max_{\xi \in \cup_{i = 1}^m X_j} \| g_j(\xi) \|_2, $ which is well-defined due to Lemma~\ref{lemma:bounded_subgradient}.
	\noindent Using the definition of $ d_i(k) $ -- see Algorithm~\ref{algo:distributed_algorithm} -- into the second term in the right-hand side of~\eqref{eq:ineq_1}, we then have that (via double stochasticity of $ A $)
	\small
	\begin{align}
	2c(k) &\sum_{i = 1}^m  d_i(k)^T (\bar{v}(k+1) - x^\star) \nonumber \\
	 &=  2c(k) \sum_{i = 1}^m  g_i(z_i(k))^T (\bar{v}(k+1) - x^\star). 
	\label{eq:equality_1}
	\end{align}
	\normalsize
	Moreover, by adding and subtracting $ x_i(k+1) $ and $ z_i(k) $ for all $ i = 1, \ldots, m, $ into the right-hand side of~\eqref{eq:equality_1} we obtain 
	\small
	\begin{align}
	2c(k)&\sum_{i = 1}^m g_i(z_i(k))^T (\bar{v}(k+1) - x^\star) \nonumber \\ =  &2c(k) \sum_{i = 1}^m  g_i(z_i(k))^T (\bar{v}(k+1) - x_i(k+1)) \nonumber \end{align} \begin{align} 
	+  &2c(k) \sum_{i = 1}^m g_i(z_i(k))^T (x_i(k+1) - z_i(k))\nonumber \\ + &2c(k) \sum_{i = 1}^m g_i(z_i(k))^T (z_i(k) - x^\star).
	\label{eq:equality_2}
	\end{align}
	\normalsize
	Consider now the right-hand side of~\eqref{eq:equality_2}. The left-most term can be lower-bounded as
	\small
	\begin{align}
	2c(k) &\sum_{i = 1}^m  g_i(z_i(k))^T (\bar{v}(k+1) - x_i(k+1)) \nonumber \\ &\geq -2 c(k) L \sum_{i = 1}^m \| \bar{v}(k+1))- x_i(k+1) \|_2,
	\label{eq:ineq_first_term}
	\end{align}
	\normalsize
	by the Cauchy-Schwartz inequality. As for the middle term, we have that
	\small
	\begin{align}
	&2c(k)  \sum_{i = 1}^m g_i(z_i(k))^T (x_i(k+1) - z_i(k)) \nonumber \\  &\geq -2c(k) L \sum_{i = 1}^m  \| e_i(k+1) \|_2 \nonumber \\
	& \geq -\alpha_1(k) \sum_{i = 1}^m \| e_i(k+1) \|_2^2 - m\frac{L^2}{\alpha_1(k)} c(k)^2
	\label{eq:ineq_middle_term}
	\end{align}
	\normalsize
	where the first inequality follows from the Cauchy-Schwartz inequality and the  definition $ e_i(k) $ in~\eqref{eq:error_def}. For the second inequality, we employed the relation $ 2xy \leq x^2 + y^2 $ with $ x = \frac{L}{\sqrt{\alpha_1(k)}} c(k) $ and $ y = \sqrt{\alpha_1(k)} \| e_i(k+1) \|_2  $ for some $ \alpha_1(k) \in (0,1) $, $ k \in \mathbb{N} $.
	
	Similarly, the right-most term of~\eqref{eq:equality_2} can be manipulated to yield
	\small
	\begin{align}
	2c(k) &\sum_{i = 1}^m g_i(z_i(k))^T (z_i(k) - x^\star) \nonumber\\ &\geq  2c(k) \sum_{i = 1}^m \Big( f_i(z_i(k)) - f_i(x^\star) \Big) \nonumber \\ &= 2c(k) \sum_{i = 1}^m \Big(  f_i(z_i(k)) - f_i(\bar{v}(k+1)) \Big) \nonumber \\ &+ 2c(k)\sum_{i = 1}^m \Big( f_i(\bar{v}(k+1)) - f_i(x^\star) \Big)
	\label{eq:long_ineq_1}
	\end{align}
	\normalsize
	where the inequality follows from the definition of the subgradient for a convex function, and the equality by adding and subtracting $ f_i(\bar{v}(k+1)) $. The first term in the right-hand side of~\eqref{eq:long_ineq_1} can be lower bounded as
	\small
	\begin{align}
	&2c(k) \sum_{i = 1}^m \Big(  f_i(z_i(k)) - f_i(\bar{v}(k+1)) \Big) \nonumber \\ &\geq 
	 -2c(k)L \sum_{i = 1}^m \| z_i(k) - \bar{v}(k+1) \|_2 \nonumber \end{align} \begin{align}
	 &\geq -2c(k) L \left( \sum_{i = 1}^m \left(  \| e_i(k+1) \|_2 + \| x_i(k + 1) - \bar{v}(k+1) \|_2 \right) \right) \nonumber \\ 
	 &\geq -\alpha_2(k) \sum_{i = 1}^m \| e_i(k+1) \|_2^2 - m \frac{L^2}{\alpha_2(k)} c(k)^2 \nonumber \\ & -2c(k) L \sum_{i = 1}^m \| x_i(k + 1) - \bar{v}(k+1) \|_2
	 \label{eq:long_ineq_2}
	\end{align}
	\normalsize
	where the first inequality follows from the relation $ x \geq -|x| $, for all $ x \in\mathbb{R} $, and from item $ iii)$ of Lemma~\ref{lemma:bounded_subgradient}, and the second inequality by adding and subtracting $ x_i(k+1) $, for all $ i =1, \ldots, m $, and then using triangle inequality. The last inequality follows from $ 2xy \leq x^2 + y^2 $ with $ x =  \frac{L}{\sqrt{\alpha_2(k)}} c(k)  $ and $ y = \sqrt{\alpha_2(k)} \| e_i(k+1)\|_2 $ for some $ \alpha_2(k) \in (0,1) $, $ k \in \mathbb{N} $. Substituting~\eqref{eq:long_ineq_2} into~\eqref{eq:long_ineq_1}
	\small
	\begin{align}
	2c(k)& \sum_{i = 1}^m g_i(z_i(k))^T (z_i(k) - x^\star) \nonumber \\  &\geq  -\alpha_2(k) \sum_{i = 1}^m \| e_i(k+1) \|_2^2 - m \frac{L^2}{\alpha_2(k)} c(k)^2 \nonumber \\
	& -2c(k) L \sum_{i = 1}^m \| x_i(k + 1) - \bar{v}(k+1) \|_2 \nonumber \\  & +   2c(k) \sum_{i = 1}^m \Big(  f_i(\bar{v}(k+1)) - f_i(x^\star)\Big). 
	\label{eq:ineq_last_term}
	\end{align}
	\normalsize
	 Substituting~\eqref{eq:ineq_1},~\eqref{eq:ineq_2},~\eqref{eq:ineq_first_term},~\eqref{eq:ineq_middle_term},~\eqref{eq:ineq_last_term} into~\eqref{eq:main_ineq} 
	 \small
	\begin{align}
	2 c(k) &\sum_{i = 1}^m (f_i(\bar{v}(k+1)) - f_i(x^\star)) + \sum_{i = 1}^m \| x_i(k+1) - x^\star \|_2^2 \nonumber \\  + &\Big( 1-\alpha_1(k) -\alpha_2(k) \Big) \sum_{i = 1}^m  \| e_i(k+1) \|_2^2  \nonumber\\
		\leq &\sum_{i = 1}^m \| x_i(k)- x^\star \|_2^2  + m L^2  \Big( \frac{\alpha_1(k) + \alpha_2(k)}{\alpha_1(k) \alpha_2(k)} \Big) c(k)^2  \nonumber  \\
	&+ 6 c(k) L \sum_{i = 1}^m \| x_i(k+1) - \bar{v}(k+1) \|_2.
	\label{eq:ineq_manupualtion_wthout_k}
	\end{align}
	\normalsize
	Summing~\eqref{eq:ineq_manupualtion_wthout_k} from $ k = 0 $ to $ k = N $, and using Lemma~\ref{lemma:main_ineq}, item $ iii) $, with $ \bar{L} = 3L $, the desired inequality~\eqref{eq:main_ineq_with_k} follows. This concludes the proof of item $ ii) $. 
\end{pf}

 Note that for any $ \beta_1 \in (0,1) $, the sequences $ (\alpha_1(k))_{k \in \mathbb{N}} $ and $ (\alpha_2(k))_{k \in \mathbb{N}} $ can be chosen to guarantee that $ 1-\alpha_1(k)-\alpha_2(k)-\beta_1 \geq 0 $ for all $ k \in \mathbb{N} $. For instance, one particular choice is $ \alpha_1(k) = \alpha_2(k) = \alpha $ with $ 1-\beta_1-2\alpha > 0 $. Three immediate consequences of Lemma~\ref{lemma:main_ineq} are presented in the following proposition. 

\begin{prop}
	Consider Assumptions~\ref{assump:Convexity_Compt}--\ref{assump:time_varying_step_size_1_t}. The following statements hold
	\begin{itemize}
		\item[$ i) $] We have that $  \sum_{k = 0}^ \infty \sum_{i = 1}^m \| e_i(k) \|_2^2 < \infty $;
		\item[$ ii) $] For all $ i = 1, \ldots, m $, we have that $ \lim_{k \rightarrow \infty} \| e_i(k) \|_2 = 0 $;
		\item[$ iii) $] For all $ i = 1, \ldots, m $, $ 
		\lim_{k \rightarrow \infty} \| x_i(k)  - v(k) \|_2 = 0.$
	\end{itemize}
	\label{proposition:intermediate_res}
\end{prop}
\begin{pf}
	Item $ i) $: Consider Lemma~\ref{lemma:main_ineq}, item $ ii) $. Note that $ \sum_{k = 0}^N \sum_{i = 1}^m \| x_i(k+1)- x^\star \|_2 $ and $ \sum_{k = 0}^N \sum_{i = 1}^m \| x_i(k) - x^\star \|_2 $ form a telescopic series, so they can be replaced by $ \sum_{i = 1}^m \| x_i(N+1) - x^\star \|_2 $ and $ \sum_{i = 1}^m \| x_i(0) - x^\star \|_2 $, respectively. Let $ \beta_1  \in (0,1)$, choose  $ \alpha_1(k) = \alpha_2(k) = \alpha $ so that $ 1 - 2 \alpha - \beta_1  > 0 $. Observe that  $ \sum_{i = 1}^m (f_i(\bar{v}(k+1)) - f_i(x^\star)) \geq 0 $ for all $ k \in \mathbb{N} $, due to optimality of $ x^\star $, so this term can be dropped from~\eqref{eq:main_ineq_with_k}. Besides, we can also drop the term $ \sum_{i = 1}^m \| x_i(N+1) - x^\star \|_2^2 \geq 0 $ since it is non-negative and appears in the left-hand side of~\eqref{eq:main_ineq_with_k}. This yields
	\small
	\begin{equation*}
	\begin{aligned}
	&(1-2\alpha - \beta_1) \sum_{k = 0}^N \sum_{i = 1}^m  \| e_i(k+1) \|_2^2 \leq \sum_{i = 1}^m \| x_i(0)- x^\star \|_2^2 \\ & \hspace{1.5cm}+  \Bigg( m L^2 \frac{2}{\alpha}  + \beta_2 \Bigg) \sum_{k = 0}^N c(k)^2 + \beta_3.
	\end{aligned}
	\end{equation*}
	\normalsize
	\noindent Letting $ N \rightarrow \infty $, we conclude that $ \sum_{k = 0}^\infty \sum_{i = 1}^m \| e_i(k)\|_2^2 $ is finite since the sequence $ (c(k))_{k \in \mathbb{N}} $ is square-summable under Assumption~\ref{assump:time_varying_step_size_1_t} and the feasible set is compact. This concludes the proof of item $ i) $.
	
	Item $ ii) $: Follows directly from item $ i) $.
	
	Item $ iii) $: This proof follows directly from the arguments presented in~\cite[Proposition 3]{MFGP:18}, and is omitted for brevity.
\end{pf}

\subsection{Proof of Theorem~\ref{theo:convergence_algo_1}}
\label{app:Proof_theo_1}

We are now in a position to prove Theorem~\ref{theo:convergence_algo_1}. To this end, we use the inequality~\eqref{eq:ineq_manupualtion_wthout_k} and leverage on Lemma 3.4 in~\cite{BT:96} to establish convergence of the sequences $ (\| x_i(k) - x^\star \|_2)_{k \in \mathbb{N}} $,  $ i = 1, \ldots, m $, to zero for some minimiser $ x^\star $ of~\eqref{eq:main_pro}. We first present Lemma 3.4 in~\cite{BT:96}.

\begin{lemma}[\cite{BT:96}]
	Consider non-negative scalar sequences $ (\ell(k))_{k \in \mathbb{N}} $,  $ (u(k))_{k \in \mathbb{N}}  $ and $ (\zeta(k))_{k \in \mathbb{N}}  $ that satisfy the recursion $
	\ell(k+1) \leq  \ell(k) - u(k) + \zeta(k).
$
	If $ \sum_{k = 0}^\infty \zeta(k) < \infty $, then the sequence $ (\ell(k))_{k \in \mathbb{N}} $ converges and the sequence $ (u(k))_{k \in \mathbb{N}} $ is summable.
	\label{lemma:supermartingale}
\end{lemma}

Consider inequality~\eqref{eq:ineq_manupualtion_wthout_k}, and choose $ \alpha_1(k), \alpha_2(k) $ and $ \beta_1 $ as in the proof of Proposition~\ref{proposition:intermediate_res} item $ i) $. We now drop the term involving $ (1-2\alpha) \sum_{i = 1}^m \| e_i(k+1)\|_2^2 $ as it appears on the left-hand side of the inequality and is non-negative so that we obtain 
\small
\begin{align}
 \sum_{i = 1}^m \| &x_i(k+1) - x^\star \|_2^2 \nonumber \leq  \sum_{i = 1}^m \| x_i(k)- x^\star \|_2^2 \nonumber \\ &- 2 c(k) \sum_{i = 1}^m (f_i(\bar{v}(k+1)) - f_i(x^\star)) +  \frac{2m L^2 }{\alpha} c(k)^2  \nonumber \\
&+ 6 c(k) L \sum_{i = 1}^m \| x_i(k+1) - \bar{v}(k+1) \|_2.
\label{eq:pre_martingale}
\end{align}
\normalsize
With reference to Lemma~\ref{lemma:supermartingale} and considering inequality~\eqref{eq:pre_martingale}, we set $\ell(k) = \sum_{i = 1}^m \| x_i(k) - x^\star \|_2^2,$ and 
\small
\begin{align}
&\zeta(k) = \frac{2  m L^2}{\alpha} c(k)^2 + 6c(k) L \sum_{i = 1}^m \| x_i(k+1) - \bar{v}(k+1) \|_2, \nonumber \\
&u(k) = 2 c(k) \big( f(\bar{v}(k+1)) - f(x^\star) \big).
\label{eq:seq_martingales}
\end{align}
\normalsize

By Lemma~\ref{lemma:literature_results}, item iii), with $\bar{L} = 3 L$, and by Proposition~\ref{proposition:intermediate_res}, item i), it follows that $ 6 L \sum_{k=1}^{\infty} c(k)  \sum_{i = 1}^m \| x_i(k+1) - \bar{v}(k+1) \| < \infty,$
hence, $\sum_{k=1}^{\infty} \zeta(k) < \infty$, as $c(k)$ is square-summable due to Assumption~\ref{assump:time_varying_step_size_1_t}, which implies that the assumptions of Lemma~\ref{lemma:supermartingale} hold. 

Therefore, we have that the sequence $ (\sum_{i = 1}^m \| x_i(k) - x^\star \|_2^2)_{k \in \mathbb{N}} $ converges, which implies that  $ (\sum_i \| x_i(k) - x^\star \|_2)_{k \in \mathbb{N}} $ also converges. To see this, note that, by continuity of the square-root function,  $ (\sum_{i = 1}^m \| x_i(k) - x^\star \|_2^2)_{k \in \mathbb{N}} $ being a convergent sequence implies that  $ (\| X(k) - x^\star\otimes\mathbf{1}^T  \|_F)_{k \in \mathbb{N}} $ also converges, where, for a fixed $ k \in \mathbb{N} $, $ X(k) $ is a $ n \times m $ matrix whose $ i$-th column is given by $ x_i(k) $, and $ \otimes $ represents the Kronecker product. Moreover, note that the set of $ n \times m $ matrices can be equipped with the norm $ \sum_{i = 1}^m \| x_i \|_2 $, where $ x_i $, $ i = 1, \ldots, m $, is the $ i$-th column of a generic element $ X \in \mathbb{R}^{n \times m} $. Since all norms in finite-dimensional spaces are equivalent, we conclude that the sequence $ (\sum_{i = 1}^m \| x_i(k) - x^\star \|_2)_{k \in \mathbb{N}} $ also converges. An alternative but more tedious justification of this argument can be found in \cite{MFGP:18}. 

%
By Lemma~\ref{lemma:supermartingale}, we also have that
$
\sum_{k = 1}^\infty c(k) \big( f(\bar{v}(k+1)) - f(x^\star) \big) < \infty.
$The latter implies that $\liminf_{k \rightarrow \infty} (f(\bar{v}(k+1)) - f(x^\star)) = 0.$ Therefore, there exists a subsequence of $ (f(\bar{v}(k+1)) - f(x^\star))_{k \in \mathbb{N}} $ that converges to zero. Since the function $ f(x) $ is continuous (by convexity) there exists some minimizer $ x^\star $ such that a subsequence of $ (\| \bar{v}(k) - x^\star \|_2)_{k \in \mathbb{N}} $ converges to zero. Moreover, we obtain $
\sum_{i = 1}^m \| x_i(k) -  x^\star \|_2 \leq \sum_{i = 1}^m \| \bar{v}(k) - x^\star \|_2 + \mu \sum_{i=1}^{m} \| x_i(k) - v(k) \|_2.
$
by adding and subtracting $ \bar{v}(k) $, then applying triangle inequality and invoking Lemma~\ref{lemma:literature_results}, item $ i) $. 

Note that $ (\| \bar{v}(k) - x^\star \|_2)_{k \in \mathbb{N}} $ converges to zero across a subsequence and $ (\sum_{i = 1}^m \| x_i(k) - v(k) \|_2)_{k \in \mathbb{N}} $ converges to zero (due to Proposition~\ref{proposition:intermediate_res}, item $ iii) $) hence we can find a subsequence of $ (\sum_{i = 1}^m \| x_i(k) - x^\star\|_2)_{k \in \mathbb{N}} $ that converges to zero. However, we have shown by means of Lemma~\ref{lemma:supermartingale} that the sequence $ (\sum_{i = 1}^m \| x_i(k) - x^\star \|_2)_{k \in \mathbb{N}} $ converges; as a result it should converge to zero since every Cauchy sequence has a unique limit point. To conclude the proof, note that, for all $ k \in \mathbb{N} $ and for all $ j = 1, \ldots, m $, $ \| x_{j}(k) - x^\star \|_2 \leq  \sum_{i = 1}^m \| x_i(k) - x^\star \|_2 $, so we conclude that the sequences $ (\| x_{j}(k) - x^\star \|_2)_{k \in \mathbb{N}} $, $ j = 1, \ldots, m  $, converge to zero. This concludes the proof.

\subsection{Proof of Theorem~\ref{theo:convergence_rate}}
\label{app:Proof_theo_2}
Consider Assumption~\ref{assump:time_varying_step_size_1_sqrt}. We drop the constant $ \eta $ for simplicity of exposition, but general choices $ \frac{\eta}{\sqrt{k+1}} $, $ \eta > 0 $, are also applicable. Let $ (\hat{v}(k))_{k \in \mathbb{N}} $ be the running average sequence associated with $ (\bar{v}(k))_{k \in \mathbb{N}} $ (definition is analogous to $ (\hat{x}_i(k))_{k \in \mathbb{N}} $ in~\eqref{eq:running_average}). Note that since $ \cap_{i = 1}^m X_i $ is assumed to be convex, we have that $ \hat{v}(k) $ is feasible for all $ k \in \mathbb{N} $ (see also the discussion below \eqref{eq:def_v_bar}). We have that  
\small
\begin{align}
&\left| \sum_{i=1}^m f_i(\hat{x}_i(k + 1)) - f(x^\star) \right|
\leq f(\hat{v}(k + 1)) - f(x^\star) \nonumber \\ &+ L \sum_{i = 1}^m \| \hat{x}_i(k+1) - \hat{v}(k+1) \|_2,
\label{eq:first_ineq_convergence_rate_proof}
\end{align}
\normalsize
which follows from triangle inequality and Lemma~\ref{lemma:bounded_subgradient}, item $ iii) $. Note that the first term in the right-hand side of~\eqref{eq:first_ineq_convergence_rate_proof} does not involve an absolute value due to feasibility of the sequence $ (\hat{v}(k))_{k \in \mathbb{N}} $, which in turn implies that $ f(\hat{v}(k+1)) \geq f(x^\star) $.

To facilitate subsequent statements, we change the notation in Lemma~\ref{lemma:main_ineq}, item $ ii) $, by replacing $ k $ by $ r $, and $ N $ by $ k $. The inequality with this modified notation is repeated here for clarity. Indeed, we have that for all $ k \in \mathbb{N} $
\small
\begin{align}
&2 \sum_{r = 0}^k c(r) \sum_{i = 1}^m (f_i(\bar{v}(r+1)) - f_i(x^\star)) \nonumber \\  &+  \sum_{r = 0}^k (1-\alpha_1(r) - \alpha_2(r) - \beta_1) \sum_{i = 1}^m  \| e_i(r+1) \|_2^2  \nonumber  \\  & + \sum_{r = 0}^k \sum_{i = 1}^m \| x_i(r+1) - x^\star \|_2^2    \leq \sum_{r = 0}^k \sum_{i = 1}^m \| x_i(r)- x^\star \|_2^2 \nonumber \\ &+  \sum_{r = 0}^k \Big( m L^2  \frac{\alpha_1(r) + \alpha_2(r)}{\alpha_1(r) \alpha_2(r)} + \beta_2 \Big) c(r)^2  + \beta_3,
\label{eq:main_ineq_with_k_mod}
\end{align}
\normalsize
where $ (\alpha_1(r))_{r \in \mathbb{N}} $ and $ (\alpha_2(r))_{r \in \mathbb{N}} $ are sequences such that $ 1-\beta_1-\alpha_1(r) - \alpha_2(r) \geq 0$ for all $ r \in \mathbb{N} $.

The proofs of items $ i) $, $ ii) $ and $ iii) $ of Theorem~\ref{theo:convergence_rate} are intertwined and will be composed into two parts: we first assume that there exist constants $ d_1, d_2, d_3, d_4 > 0 $ such that~\eqref{eq:term_A_conv_rate_proof} and~\eqref{eq:term_B_conv_rate_proof} bellow are satisfied, and on this basis prove the claims of the theorem; we then return to~\eqref{eq:term_A_conv_rate_proof} and~\eqref{eq:term_B_conv_rate_proof}, and prove the existence of such constants. To this end, consider
\small
\begin{align}
&f(\hat{v}(k+1)) - f(x^\star) \leq  d_1 \frac{1}{S(k +1)}  + d_2 \frac{\sum_{r = 0}^{k} c(r)^2}{S(k + 1)} 
\label{eq:term_A_conv_rate_proof}\\
&L \sum_{i = 1}^m \| \hat{x}_i(k+1) - \hat{v}(k+1) \|_2 \leq     \frac{d_3}{S(k+1)} + d_4 \frac{\sum_{r = 0}^{k} c(r)^2}{S(k+1)}. \label{eq:term_B_conv_rate_proof} 
\end{align}
\normalsize
Note that $ S(k + 1) $ can be lower-bounded as
\small
\begin{align}
S(&k + 1) = \sum_{r = 1}^{k + 1} \frac{1}{\sqrt{r+1}} \geq \int_{2}^{k+3} \frac{1}{\sqrt{x}} dx \nonumber \\ &= 2(\sqrt{k+3} - \sqrt{2}) \geq \nu \sqrt{k + 3} \geq \nu \sqrt{k + 1} ,
\label{eq:SN_bound}
\end{align}
\normalsize
with $ \nu = 2 - \sqrt{2}$, and where we employed monotonicity of $ \frac{\sqrt{x+3}-\sqrt{2}}{\sqrt{x+1}} $ for $ x \geq 1 $. Moreover, we have that
\small
\begin{align}
\sum_{r = 0}^{k} c(r)^2 = &\sum_{r =0}^{k} \frac{1}{r + 1} = \sum_{r = 1}^{k+1} \frac{1}{r}   \nonumber \\ &\leq  \int_{1}^{k+1} \frac{1}{x} dx + 1\leq \ln(k+1) + 1.
\label{eq:c_k_over_S_k}
\end{align}
\normalsize
The result of the Theorem~\ref{theo:convergence_rate}, item $ iii) $, follows then from~\eqref{eq:first_ineq_convergence_rate_proof} by substituting~\eqref{eq:term_A_conv_rate_proof}--\eqref{eq:c_k_over_S_k}, and setting $ B_1 = \sum_{i = 1}^4 \frac{d_i}{\nu} $ and $ B_2 = \frac{d_2}{\nu} + \frac{d_4}{\nu} $. 
Since \eqref{eq:term_B_conv_rate_proof} is valid for all $ i = 1, \ldots, m $,  we have that (via a direct application of triangle inequality)
$
\| \hat{x}_i(k) - \hat{x}_j(k) \|_2 \leq  \sum_{i = 1}^m \| \hat{x}_i(k) - \hat{v}(k) \| +  \sum_{i = 1}^m \| \hat{x}_j(k) - \hat{v}(k) \|, 
$
which due to \eqref{eq:SN_bound} and \eqref{eq:c_k_over_S_k} then implies that the sequence $  (\| \hat{x}_i(k) - \hat{x}_j(k) \|_2)_{k \in \mathbb{N}}$ converges to zero at a rate $ \mathcal{O}(\frac{\ln k}{\sqrt{k}}) $. This concludes the proof of item $ i) $. 

Moreover, these relations also imply that the set of accumulation points of the sequence $ (\hat{v}(k))_{k \in \mathbb{N}} $ coincides to that of the sequences $ (\hat{x}_i(k))_{k \in \mathbb{N}} $, $ i  = 1, \dots, m $. Hence, we conclude that all accumulation points of $ (\hat{x}_i(k))_{k \in \mathbb{N}} $ are feasible due to the fact that all accumulation points of $ (\hat{v}(k))_{k \in \mathbb{N}} $ are in $ \cap_{i = 1}^m X_i $ and the latter is a closed set, thus concluding the proof of item $ ii) $. This concludes the proof of Theorem~\ref{theo:convergence_rate}. 

\subsection*{Derivation of~\eqref{eq:term_A_conv_rate_proof}}

We first construct an upper-bound for the term on the left-hand side of~\eqref{eq:term_A_conv_rate_proof}. In fact, observe that \small
\begin{align}
f(&\hat{v}(k+1)) - f(x^\star) = f\left(\frac{1}{S(k+1)}\sum_{r = 1}^{k+1} c(r) \bar{v}(r)\right) - f(x^\star) \nonumber  \\ &\leq  \sum_{r = 1}^{k+1} \frac{c(r)}{S(k+1)} f(\bar{v}(r)) - f(x^\star) \nonumber \\ 
&= \sum_{r = 0}^k \frac{c(r+1)}{S(k+1)} \sum_{i = 1}^m (f_i(\bar{v}(r+1)) - f_i(x^\star)) \nonumber  \\
&\leq \sum_{r = 0}^k \frac{c(r)}{S(k+1)} \sum_{i = 1}^m (f_i(\bar{v}(r+1)) - f_i(x^\star)),
\label{eq:manipulation_LHS_term_A}
\end{align}
\normalsize
where the first equality follows by definition of $ \hat{v}(k+1) $, the first inequality by convexity of $ f $, the second equality by using the fact that $ f = \sum_{i = 1}^m f_i $ and changing the summation index, and the second inequality by using the fact that $ c(r+1) = \frac{1}{\sqrt{r+1}} \leq  \frac{1}{\sqrt{r}} = c(r)   $ for all $ r \in \mathbb{N} $.

In light of~\eqref{eq:main_ineq_with_k_mod},  for any $ \beta_1 \in (0,1) $, a valid choice for the sequences $ (\alpha_1(k))_{k \in \mathbb{N}} $ and $ (\alpha_2(k))_{k \in \mathbb{N}} $ is $ \alpha_1(k) = \alpha_2(k) = \alpha(k) $, where $ \alpha(k) =  a \Big( 1 - \frac{1}{\sqrt{k+1}} \Big)$; to ensure that $ 1-\beta_1 - \alpha_1(k) - \alpha_2(k) \geq 0 $ as required by Lemma~\ref{lemma:main_ineq}, item $ ii) $, it suffices to set $ a = (1-\beta_1)/2 $. Under these choices we have that
\small
\begin{equation}
1-\beta_1 - 2\alpha(k) = \frac{1-\beta_1}{\sqrt{k+1}} = (1-\beta_1) c(k).
\label{eq:alpha_error_ineq}
\end{equation}
\normalsize
Consider now~\eqref{eq:main_ineq_with_k_mod} with the above choices for $ \alpha_1(k) $ and $ \alpha_2(k) $. Note that the series $ \sum_{r = 0}^{k} \sum_{i=1}^m \| x_i(r+1) - x^\star \|_2  $ and $ \sum_{r = 0}^{k} \sum_{i=1}^m \| x_i(r) - x^\star \|_2  $ are telescopic, thus all intermediate terms cancel. We now drop the terms involving $ \| e_i(r+1) \|_2^2 $ and $ \| x_i(k+1) - x^\star \|_2 $ as they are non-negative, and then divide the resulting expression by $ 2 S(k + 1) = 2 \sum_{r = 1}^{k + 1} \frac{1}{\sqrt{r+1}} $ to obtain the following upper bound on the right-hand side of~\eqref{eq:manipulation_LHS_term_A}
\small
\begin{align}
&\sum_{r = 0}^{k} \frac{c(r)}{S(k + 1)} \sum_{i = 1}^m (f_i(\bar{v}(r+1)) - f_i(x^\star)) \nonumber \\ &\leq \frac{\sum_{i = 1}^m \| x_i(0)- x^\star \|_2^2}{2S(k + 1)} +  \frac{\beta_3}{2S(k + 1)} \nonumber\\
&+ \frac{\beta_2}{2} \sum_{r = 0}^{k} \frac{c(r)^2}{S(k + 1)}    + m L^2 \frac{1}{S(k+1)} \sum_{r = 0}^{k}  \frac{ c(r)^2}{\alpha(r)}.
\label{eq:ineq_crucial_lemma_rate_convergence}
\end{align}
\normalsize
By the right-hand side of~\eqref{eq:ineq_crucial_lemma_rate_convergence}, we obtain~\eqref{eq:term_A_conv_rate_proof} with $ 
d_1 = \frac{4mD^2 + \beta_3}{2}, \quad d_2 = \frac{\beta_2}{2} + \frac{4mL^2}{a}.
$
where, by Assumption~\ref{assump:Convexity_Compt},  $ \sum_{i = 1}^m \| x_i(0) - x^\star \|_2^2 \leq 4mD^2 $, with $ D $ defined as in Lemma~\ref{lemma:literature_results}, item $ i) $. Moreover, we used the fact that $
\frac{c(r)^2}{\alpha(r)} = \frac{1}{a} \frac{\sqrt{r+1}}{\sqrt{r+1} - 1} \frac{1}{r + 1} \leq \frac{4}{a} c(r)^2, 
$
due to monotonicity of $ \frac{\sqrt{x+1}}{\sqrt{x+1} - 1} $. 
\subsection*{Derivation of~\eqref{eq:term_B_conv_rate_proof}}
Similarly to the derivation of~\eqref{eq:term_A_conv_rate_proof}, we apply the definition of both $ \hat{x}_i(k) $, $ i = 1, \ldots, m $, and $ \hat{v}(k) $ to upper-bound the left-hand side of~\eqref{eq:term_B_conv_rate_proof} as \small
\begin{align}
L &\sum_{i = 1}^m \| \hat{x}_i(k+1) - \hat{v}(k+1) \|_2 \nonumber \\ &= L \sum_{i = 1}^m \left\| \frac{1}{S(k+1)} \sum_{r = 1}^{k+1} c(r)\Big( x_i(r) - \bar{v}(r)  \Big) \right\|_2 \nonumber \\ &\leq \frac{L\mu}{S(k+1)} \sum_{r = 1}^{k+1} c(r) \sum_{i = 1}^m \| x_i(r) - v(r) \|_2,
\label{eq:first_term_term_B_theo_2}
\end{align}
\normalsize
where the inequality follows from convexity of the norm. We will now construct an upper-bound on the right-hand side of~\eqref{eq:first_term_term_B_theo_2}. To this end, note that
\small
\begin{align}
&\frac{L\mu}{S(k+1)} \sum_{r = 1}^{k+1} c(r) \sum_{i = 1}^m \| x_i(r) - v(r)\|_2 \nonumber \\ &= \frac{L \mu c(1)}{S(k+1)} \sum_{i = 1}^m \| x_i(1) - v(1) \|_2\nonumber \\ &+   \frac{L \mu}{S(k+1)} \sum_{r = 2}^{k+1} c(r) \sum_{i = 1}^m \| x_i(r) - v(r)\|_2.
\label{eq:ineq_intermidiate_conv_rate_proof}
\end{align}
\normalsize
 We now invoke Lemma~\ref{lemma:literature_results}, item $ ii) $ -- with $ r $ in the place of $ k $, and $ t  $ in the place of $ r $ -- for the last term on the right-hand side of~\eqref{eq:ineq_intermidiate_conv_rate_proof} so that
 \small
\begin{align}
&\sum_{r = 2}^{k+1} c(r) \sum_{i=1}^m \| x_i(r) - v(r) \|_2 \nonumber\\ &= \sum_{r = 1}^{k} c(r +1) \sum_{i=1}^m \| x_i(r+1) - v(r+1) \|_2 \nonumber \\
&\leq  2 \sum_{r = 0}^{k} c(r) \sum_{i = 1}^{m} \| e_i(r+1) \|_2  + m\lambda \sum_{i=1}^m \| x_i(0) \|_2 \sum_{r = 0}^{k} c(r) q^{r} \nonumber \\
&+m\lambda\sum_{r = 1}^{k} c(r+1) \sum_{t = 0}^{r-1} q^{r-t-1} \sum_{i = 1}^{m} \| e_i(t+1) \|_2
\label{eq:ineq_intermidiate_conv_rate_proof_2}
\end{align}
\normalsize
where we added the term corresponding to $ r = 0 $ and used the fact that $ c(r+1) \leq c(r) $ for all $ r \in \mathbb{N} $, in first two terms on the right-hand side of~\eqref{eq:ineq_intermidiate_conv_rate_proof_2}.  We analyse each term in the right-hand side of~\eqref{eq:ineq_intermidiate_conv_rate_proof_2} separately. First, observe that 
\small
\begin{align}
2 \sum_{r = 0}^k c(r) \sum_{i = 1}^m \| e_i(r+1) &\|_2 \leq \sum_{r = 0}^k c(r)^2 + \sum_{i = 1}^m \| e_i(r+1) \|_2^2,
\label{eq:first_term_B_theo_2_final}
\end{align}
\normalsize
using the identity $ 2xy \leq x^2 + y^2 $. The intermediate term in the right-hand side of~\eqref{eq:ineq_intermidiate_conv_rate_proof_2} can be manipulated to yield
\small
\begin{align}
m \lambda \sum_{i = 1}^m \| x_i(0) \|_2 \sum_{r = 0}^k c(r) q^r \leq \frac{m^2 \lambda D}{1-q}, 
\label{eq:last_term_3}
\end{align}
\normalsize
since $ c(r) \leq 1 $ for all $ r \in \mathbb{N}\cup \{ 0\} $, $ \| x_i(0) \|_2 \leq D $ (Lemma~\ref{lemma:bounded_subgradient}) for all $ i = 1, \ldots, m $, and using the closed-form expression for the sum of geometric series as $ q \in (0,1) $. We deal with the last term in~\eqref{eq:ineq_intermidiate_conv_rate_proof_2} in several steps. We start by expanding the terms to obtain
\small
\begin{align}
&\sum_{r = 1}^k c(r+1) \sum_{t = 0}^{r-1} q^{r-t-1} \sum_{i = 1}^m \| e_i(t+1) \|_2 \nonumber \\ &= c(2)  \sum_{i = 1}^m \|e_i(1)\|_2 
+c(3) \left( q \sum_{i = 1}^m \| e_i(1) \|_2 \nonumber \sum_{i = 1}^m \| e_i(2) \|_2\right)  \\ & + \ldots + c(k+1) \left( \sum_{t = 1}^{k}q^{k-t} \sum_{i = 1}^m \| e_i(t) \|_2 \right). \label{eq:last_term_1}
\end{align}
\normalsize
We now collect the terms containing the error vector $ e_i(r) $, $ r = 1, \ldots, k $, to obtain
\small
\begin{align}
&m\lambda\sum_{r = 1}^{k} c(r+1) \sum_{t = 0}^{r-1} q^{r-t-1} \sum_{i = 1}^{m} \| e_i(t+1) \|_2 \nonumber  \\
&= m \lambda\sum_{i = 1}^m \| e_i(1) \|_2 \Big(c(2) + q c(3) + \ldots \nonumber  \\ &+ q^{k-1} c(k+1) \Big) + \ldots  + \sum_{i = 1}^m \| e_i(k) \|_2 c(k+1) \nonumber \end{align} \begin{align}
&\leq \frac{m\lambda}{1-q} \sum_{r = 1}^{k} c(r+1) \sum_{i = 1}^m \| e_i(r) \|_2\nonumber\\ &\leq \frac{m \lambda}{1-q} \sum_{r = 1}^k c(r) \sum_{i = 1}^m \| e_i(r) \|_2 \leq   \frac{m \lambda}{2(1-q)} \sum_{r = 0}^k c(r)^2 \nonumber \\ & \quad  \quad \quad\quad+  \frac{m \lambda}{2(1-q)} \sum_{r = 0}^k \sum_{i = 1}^m \| e_i(r+1) \|_2^2  
 \label{eq:last_term_2}
\end{align}
\normalsize
where in the first inequality we used the fact that $ q \leq \frac{1}{1-q}$ and $ 1 \leq \frac{1}{1-q} $ for any $ q \in (0,1) $, while in the second inequality we used the fact that $ c(r + 1) \leq c(r) $. To obtain the last inequality we applied the relation $ 2xy \leq x^2 + y^2 $ with $ x = c(r) $ and $ y = \| e_i(r+1) \|_2 $, and then added the non-negative terms involving $ c(0)^2  $ and $ \sum_{i = 1}^m \| e_i(k+1) \|_2^2 $. Substituting~\eqref{eq:ineq_intermidiate_conv_rate_proof}--\eqref{eq:last_term_3} and~\eqref{eq:last_term_2} into~\eqref{eq:first_term_term_B_theo_2} we have that
\small
\begin{align}
&L \sum_{i = 1}^m \| \hat{x}_i(k+1) - \hat{v}(k+1) \|_2 \nonumber \\ &\leq L \mu \left( 1 + \frac{m \lambda}{2(1-q)} \right) \frac{\sum_{r = 0}^k c(r)^2}{S(k+1)} \nonumber \\ &+ \Bigg(   m \lambda + 2 c(1) \Bigg)     \frac{L\mu m D}{S(k+1)} \nonumber \\ &
+  \frac{L\mu}{S(k+1)} \left( 1 + \frac{m \lambda}{2(1-q)} \right)   \sum_{r = 0}^k \sum_{i = 1}^m \| e_i(r+1) \|_2^2.   
\label{eq:last_term_B} 
\end{align}
\normalsize
To obtain the result, we need to manipulate the last term in the right-hand side of~\eqref{eq:last_term_B}. To this end, we invoke~\eqref{eq:main_ineq_with_k_mod} with the same $ \beta_1  $ as in~\eqref{eq:alpha_error_ineq}, but with $ (\alpha_1(k))_{k \in \mathbb{N}} $ and $ (\alpha_2(k))_{k \in \mathbb{N}} $ such that $ \alpha_1(k) = \alpha_2(k) = \alpha $, for all $ k \in \mathbb{N} $, following the same rationale as in Proposition~\ref{proposition:intermediate_res} to obtain
\small
\begin{align}
\sum_{r = 0}^{k} &\sum_{i=1}^m \| e_i(r+1) \|_2^2 \leq \frac{\sum_{i = 1}^m \| x(0)-x^\star\|_2^2  + \beta_3}{1-\beta_1 - 2\alpha} 
\nonumber\\
 &+ \frac{1}{1-\beta_1-2\alpha} \Big( mL^2 \frac{2}{\alpha} + \beta_2 \Big) \sum_{r = 0}^{k} c(r)^2 \nonumber\\&\leq   \frac{4mD^2  + \beta_3}{1-\beta_1 - 2\alpha} \nonumber  \\
&+ \frac{1}{1-\beta_1-2\alpha} \Big( mL^2\frac{2}{\alpha} + \beta_2 \Big) \sum_{r = 0}^{k} c(r)^2.
\label{eq:conv_rate_final_3}
\end{align}
\normalsize
Substituting~\eqref{eq:conv_rate_final_3} into~\eqref{eq:last_term_B} we obtain~\eqref{eq:term_B_conv_rate_proof} with constants
\small
\begin{align*} 
d_3 &=  L\mu \Biggl[  \Bigg( 1 + \frac{m \lambda}{2(1-q)} \Bigg) \frac{4mD^2 + \beta_3}{1-\beta_1 - 2\alpha}  + m D \Bigg(  m \lambda + 2 c(1)  \Bigg) \Biggr],  \nonumber  \\
d_4 &= L \mu \Bigg( 1 + \frac{m \lambda}{2(1-q)} \Bigg) \Bigg( 1 +  \frac{1}{1-\beta_1-2\alpha} \Big( mL^2 \frac{2}{\alpha} + \beta_2 \Big) \Bigg),
\end{align*}
\normalsize
thus concluding the proof of Theorem~\ref{theo:convergence_rate}. 	
%


\begin{thebibliography}{10}
	\expandafter\ifx\csname url\endcsname\relax
	\def\url#1{\texttt{#1}}\fi
	\expandafter\ifx\csname urlprefix\endcsname\relax\def\urlprefix{URL }\fi
	\expandafter\ifx\csname href\endcsname\relax
	\def\href#1#2{#2} \def\path#1{#1}\fi
	
	\bibitem{MG:12}
	G.~Mateos, G.~B. Giannakis, {Distributed recursive least-squares: Stability and
		performance analysis}, IEEE Transactions on Signal Processing 60~(7) (2012)
	3740--3754.
	
	\bibitem{BMG:14}
	B.~Baingana, G.~Mateos, G.~B. Giannakis, {Proximal-Gradient Algorithms for
		Tracking Cascades Over Social Networks}, IEEE Journal of Selected Topics in
	Signal Processing 8~(4) (2014) 563--575.
	
	\bibitem{MBCF:07}
	S.~Martinez, F.~Bullo, J.~Cortes, E.~Frazzoli, {On Synchronous Robotic Networks
		- Part I: Models, Tasks and Complexity}, IEEE Transactions on Automatic
	Control 52~(12) (2007) 2199--2213.
	
	\bibitem{BCCZ:15}
	S.~Bolognani, R.~Carli, G.~Cavraro, S.~Zampieri, {Distributed Reactive Power
		Feedback Control for Voltage Regulation and Loss Minimization}, IEEE
	Transactions on Automatic Control 60~(4) (2015) 966--981.
	
	\bibitem{ZM:12}
	M.~Zhu, S.~Martinez, {On distributed convex optimization under inequality and
		equality constraints}, IEEE Transactions on Automatic Control 57~(1) (2012)
	151--164.
	
	\bibitem{NO:15}
	A.~Nedi{\'{c}}, A.~Olshevsky, {Distributed optimization over time-varying
		directed graphs}, IEEE Transactions on Automatic Control 60~(3) (2015)
	601--615.
	
	\bibitem{XK:17}
	C.~Xi, U.~A. Khan, {Distributed Subgradient Projection Algorithm Over Directed
		Graphs}, IEEE Transactions on Automatic Control 62~(8) (2017) 3986--3992.
	
	\bibitem{MFGP:18}
	K.~Margellos, A.~Falsone, S.~Garatti, M.~Prandini, {Distributed Constrained
		Optimization and Consensus in Uncertain Networks via Proximal Minimization},
	IEEE Transactions on Automatic Control 63~(5) (2018) 1372--1387.
	
	\bibitem{LWY:19}
	S.~Liang, L.~Wang, G.~Yin, {Distributed quasi-monotone subgradient algorithm
		for nonsmooth convex optimization over directed graphs}, Automatica 101
	(2019) 175--181.
	
	\bibitem{Bia:2016}
	P.~Bianchi, {Ergodic convergence of a stochastic proximal point algorithm},
	SIAM Journal on Optimization 26~(4) (2016) 2235--2260.
	
	\bibitem{PN:18}
	A.~Patrascu, I.~Necoara, {Nonasymptotic convergence of stochastic proximal
		point methods for constrained convex optimization}, Journal of Machine
	Learning Research 18 (2018) 1--42.
	
	\bibitem{JKJJ:08}
	B.~Johansson, T.~Keviczky, M.~Johansson, K.~H. Johansson, {Subgradient methods
		and consensus algorithms for solving convex optimization problems},
	Proceedings of the IEEE Conference on Decision and Control~(5) (2008)
	4185--4190.
	
	\bibitem{NO:09}
	A.~Nedic, A.~Ozdaglar, {Approximate primal solutions and rate analysis for dual
		subgradients methods}, SIAM Journal on Optimization 33~(5) (2009) 2295--2317.
	
	\bibitem{NOP:10}
	A.~Nedic, A.~Ozdaglar, P.~A. Parrilo, {Constrained Consensus and Optimization
		in Multi-Agent Networks}, IEEE Transactions on Automatic Control 55~(4)
	(2010) 922--938.
	
	\bibitem{LRS:16}
	P.~Lin, W.~Ren, Y.~Song,
	{Distributed
			multi-agent optimization subject to nonidentical constraints and
			communication delays}, Automatica 65 (2016) 120--131.
	
	\bibitem{LN:13}
	S.~Lee, A.~Nedi{\'{c}}, {Distributed random projection algorithm for convex
		optimization}, IEEE Journal on Selected Topics in Signal Processing 7~(2)
	(2013) 221--229.
	
	\bibitem{MA:19}
	V.~S. Mai, E.~H. Abed, {Distributed optimization over directed graphs with row
		stochasticity and constraint regularity}, Automatica 102 (2019) 94--104.
	
	\bibitem{TBA:86}
	J.~N. Tsitsiklis, D.~P. Bertsekas, M.~Athans, {Distributed Asynchronous
		Deterministic and Stochastic Gradient Optimization Algorithms}, IEEE
	Transactions on Automatic Control 31~(9) (1986) 803--812.
	
	\bibitem{BT:89}
	D.~P. Bertsekas, J.~N. Tsitsiklis, {Parallel and Distributed Computation:
		Numerical Methods}, Scientific, Athena, 1989.
	
	\bibitem{DAW:12}
	J.~C. Duchi, A.~Agarwal, M.~J. Wainwright, {Dual averaging for distributed
		optimization: Convergence analysis and network scaling}, IEEE Transactions on
	Automatic Control 57~(3) (2012) 592--606.
	
	\bibitem{RMNP:19b}
	L.~Romao, K.~Margellos, G.~Notarstefano, A.~Papachristodoulou, {Convergence
		rate analysis of a subgradient averaging algorithm for distributed
		optimisation with different constraint sets}, in: 58th Conference on Decision
	and Control, 2019, pp. 7448--7453.
	
	\bibitem{NO:09b}
	A.~Nedi{\'{c}}, A.~Ozdaglar, {Distributed Subgradient Methods for Multi-Agent
		Optimization}, IEEE Trans. Automatic Control 54~(1) (2009) 48--61.
	
	\bibitem{JMX:12}
	D.~Jakovetic, J.~M. Moura, J.~Xavier, {Distributed Nesterov-like gradient
		algorithms}, in: 51st IEEE Conference on Decision and Control, 2012, pp.
	5459--5464.
	
	\bibitem{YLY:16}
	K.~Yuan, Q.~Ling, W.~Yin, {On the convergence of decentralized gradient
		descent}, SIAM Journal on Optimization 26~(5) (2016) 1835--1854.
	
	\bibitem{LCF:16}
	J.~Lei, H.-F. Chen, H.-T. Fang, {Primal–dual algorithm for distributed
		constrained optimization}, Systems {\&} Control Letters 96 (2016) 110--117.
	
	\bibitem{SS:18}
	G.~Scutari, Y.~Sun, {Distributed Nonconvex Constrained Optimization over
		Time-Varying Digraphs}, Mathematical Programming, 2019, pp.497--544.
	
	\bibitem{TLR:12}
	K.~I. Tsianos, S.~Lawlor, M.~G. Rabbat, {Push-Sum Distributed Dual Averaging
		for convex optimization}, 2012 IEEE 51st IEEE Conference on Decision and
	Control (CDC) (2012) 5453--5458.
	
	\bibitem{LQX:17}
	S.~Liu, Z.~Qiu, L.~Xie, {Convergence rate analysis of distributed optimization
		with projected subgradient algorithm}, Automatica 83 (2017) 162--169.
		
	\bibitem{SLWY:15}
	W.~Shi, Q.~Ling, G.~Wu, W.~Yin, {EXTRA: An Exact First-Order Algorithm for
		Decentralized Consensus Optimization}, SIAM Journal on Optimization 25~(2)
	(2015) 944--966.
	
	\bibitem{QL:17}
	G.~Qu, N.~Li, {Harnessing smoothness to accelerate distributed optimization},
	IEEE Transactions on Control of Network Systems 5~(3) (2018) 1245--1260.
	
	\bibitem{VZCPS:16}
	F.~Zanella, D.~Varagnolo, A.~Cenedese, G.~Pillonetto, L.~Schenato, 
{Newton-Raphson
			Consensus for Distributed Convex Optimization}, IEEE Transactions on
	Automatic Control 61~(4) (2016) 994--1009.
	
	\bibitem{CO:12}
	A.~I. Chen, A.~Ozdaglar, {A fast distributed proximal-gradient method}, 50th Annual Allerton Conference on Communication, Control, and
	Computing, Allerton 2012, 2012, pp. 601--608.
	
	\bibitem{Ber:09}
	D.~P. Bertsekas, {Convex Optimization Theory}, Athena Scientific, 2009.
	
	\bibitem{Roc:72}
	R.~T. Rockafellar, {Convex Analysis}, Princeton University Press, 1972.
	
	\bibitem{BT:96}
	D.~P. Bertsekas, J.~Tsitsiklis, {Neuro-Dynamic Programming}, Athena Scientific,
	1996.
	
\end{thebibliography}
%

\end{document}